\documentclass[11pt]{article}

\textheight=\vsize   \textwidth=\hsize
\newif\ifproofmode
\proofmodetrue

\catcode`\@=11

\def\^#1{\if#1i{\accent"5E\i}\else{\accent"5E #1}\fi}
\def\"#1{\if#1i{\accent"7F\i}\else{\accent"7F #1}\fi}

\def\eqalign#1{\null\,\vcenter{\openup1\jot \m@th
  \ialign{\strut\hfil$\displaystyle{##}$&$\displaystyle{{}##}$\hfil
     &&\strut$\displaystyle{##}$\hfil&$\displaystyle{{}##}$
     \hfil\crcr#1\crcr}}\,}

\def\ldisplaylinesno#1{\displ@y\halign{
\hbox to\displaywidth{$\@lign\hfil\displaystyle##\hfil$}&
\kern-\displaywidth\rlap{$##$}
\tabskip\displaywidth\crcr#1\crcr}}

\def\lcoupe#1#2\label#3{\ldisplaylinesno{\refstepcounter{equation}
\label{#3}
\quad #1 \hfill & \cr
\hfill #2 \quad & \ifproofmode\llap{{\tiny\tt #3}\kern5mm}\fi
\reset@font\rm (\theequation)\cr}}

        [1994/12/02 v2.3e Standard LaTeX font definitions]
\DeclareFontFamily{OMX}{cmex}{}{}
\DeclareFontShape{OMX}{cmex}{m}{n}{
<5> <6> <7> <8> <9> 
      <10> <10.95> <12> <14.4> <17.28> <20.74> <24.88>cmex10
   }{}

\DeclareFontFamily{U}{lasy}{}
\DeclareFontShape{U}{lasy}{m}{n}{ <5> <6> <7> <8> <9> gen * lasy
      <10> <10.95> <12> <14.4> <17.28> <20.74> <24.88>lasy10  }{}
\DeclareFontShape{U}{lasy}{bx}{n}{ <-10> ssub * lasy/m/n
     <10> <10.95> <12> <14.4> <17.28> <20.74> <24.88>lasyb10  }{}
\fontencoding{U}
\DeclareRobustCommand\lyfamily
        {\not@math@alphabet\lyfamily\relax
         \fontfamily\lydefault\selectfont}     
\newcommand\lydefault{lasy}

\def\@normalsize{\@setsize\normalsize{16pt}\xiipt\@xiipt
\abovedisplayskip 12\p@ plus3\p@ minus7\p@
\belowdisplayskip \abovedisplayskip
\abovedisplayshortskip  \z@ plus3\p@
\belowdisplayshortskip  6.5\p@ plus3.5\p@ minus3\p@
\let\@listi\@listI}

\def\@ceqnnum{{\reset@font\rm (\tempequation)}}

\def\@ceqncr{{\ifnum0=`}\fi\@ifstar{\global\@eqpen\@M
    \@cyeqncr}{\global\@eqpen\interdisplaylinepenalty \@cyeqncr}}

\def\@cyeqncr{\@ifnextchar [{\@cxeqncr}{\@cxeqncr[\z@]}}

\def\@cxeqncr[#1]{\ifnum0=`{\fi}\@@ceqncr
   \noalign{\penalty\@eqpen
\vskip\jot
\vskip #1\relax}}

\def\@@ceqncr{\let\@tempa\relax
    \ifcase\@eqcnt \def\@tempa{& & &}\or \def\@tempa{& &}%
      \else \def\@tempa{&}\fi
     \@tempa \if@eqnsw\@ceqnnum\fi
     \global\@eqnswtrue\global\@eqcnt\z@\cr}

\def\ceqnarray{
\global\@eqnswtrue\m@th
\global\@eqcnt\z@\tabskip\@centering\let\\
\@ceqncr $$
\halign to\displaywidth\bgroup\@eqnsel\hskip\@centering
  $\displaystyle\tabskip\z@{{}##}$&\global\@eqcnt\@ne
  \hskip 2\arraycolsep \hfil${{}##}$\hfil
  &\global\@eqcnt\tw@ \hskip 2\arraycolsep$\displaystyle\tabskip\z@{{}##}$\hfil
   \tabskip\@centering&\llap{##}\tabskip\z@\cr}

\def
\endceqnarray{\@@ceqncr\egroup$$
\global\@ignoretrue}

\def\theequation{\arabic{section}.\arabic{equation}
}

\@addtoreset{equation}
{section}

\catcode`\@=12

\def\elle#1{L^{#1}(\Omega)}

\def\de{\delta}

\def\dys{\displaystyle}

\def\frac#1#2{{#1\over  #2}}
\def\vep{\varepsilon}
\def\r{{\bf R}}

\def\rn{{\bf R}^N}
\def\rnp{{\bf R}_{+}^N}

\def\vfi{\varphi}

\def\ude{u_2^\vep}

\newcommand{\be}{
\begin{equation}
}
\newcommand{\bel}[1]{
\begin{equation}
\label{#1}}
\newcommand{\ee}{
\end{equation}
}
\newtheorem{subn}{\name}

\newcommand{\bsn}[1]{\def\name{#1}
\begin{subn}}
\newcommand{\esn}{
\end{subn}}
\newtheorem{sub}{\name}[section]
\newcommand{\bs}{
\begin{sub}}
\newcommand{\es}{
\end{sub}}
\newcommand{\bsl}[1]{
\begin{sub}\label{#1}}
\newcommand{\bth}[1]{\def\name{Theorem}
\begin{sub}\label{t:#1}}
\newcommand{\blemma}[1]{\def\name{Lemma}
\begin{sub}\label{l:#1}}
\newcommand{\bcor}[1]{\def\name{Corollary}
\begin{sub}\label{c:#1}}
\newcommand{\bdef}[1]{\def\name{Definition}
\begin{sub}\label{d:#1}}
\newcommand{\bprop}[1]{\def\name{Proposition}
\begin{sub}\label{p:#1}}

\newcommand{\BA}{
\begin{array}}
\newcommand{\EA}{
\end{array}}
\newcommand{\BAN}{\renewcommand{\arraystretch}{1.2}
\setlength{\arraycolsep}{2pt}
\begin{array}}
\newcommand{\BAV}[2]{\renewcommand{\arraystretch}{#1}
\setlength{\arraycolsep}{#2}
\begin{array}}
\newcommand{\BSA}{
\begin{subarray}}
\newcommand{\ESA}{
\end{subarray}}
\newcommand{\BAL}{
\begin{aligned}}
\newcommand{\EAL}{
\end{aligned}}
\newcommand{\BALG}{
\begin{alignat}}
\newcommand{\EALG}{
\end{alignat}}
\newcommand{\BALGN}{
\begin{alignat*}}
\newcommand{\EALGN}{
\end{alignat*}}
\newcommand{\note}[1]{\textit{#1.}\hspace{2mm}}
\newcommand{\Proof}{\note{Proof}}

\newcommand{\qed}{\\
${}$ \hfill $\square$}

\newcommand{\abs}[1]{\left |#1\right |}

\def\angb<#1>{\langle #1 \rangle}

\newcommand{\dist}{\opname{dist}}

\newcommand{\prt}{
\partial}

\def\ga{\alpha}            
       \def\gd{\delta}

\def\Gw{\Omega}              

\newtheorem{theorem}{Theorem}[section]

\newtheorem{lemma}{Lemma}[section]

\newtheorem{remark}{Remark}[section]

\newtheorem{corollary}{Corollary}[section]

\newcommand{\rife}[1]{(\ref{#1})}

\def\al{\alpha}
\def\ga{\gamma}
\def\la{\lambda}
 \def\dist{{\rm dist}\,}
\def\be{
\begin{equation}
}
\def\ee{
\end{equation}
}
\def\dom{d_{ \Omega}(x)}
\def\qed{{\unskip\nobreak\hfil\penalty50
          \hskip2em\hbox{}\nobreak\hfil\mbox{\rule{1ex}{1ex} \qquad}
   \parfillskip=0pt
   \finalhyphendemerits=0
\par
\medskip}}

\begin{document}
 
\title {\bf  Asymptotic behaviour  for the gradient 
of large solutions to some nonlinear elliptic equations}

\author{{\bf\large Alessio Porretta}
\footnote{The author acknowledges the support of RTN european project:
FRONTS-SINGULARITIES, RTN contract: HPRN-CT-2002-00274. }
\hspace{2mm}\vspace{3mm}\\
{\it Dipartimento di Matematica,  Universit\`a di Roma Tor 
Vergata},\\
{\it Via della Ricerca Scientifica 1, 00133 Roma, Italia}\\
\vspace{3mm}\\
{\bf\large Laurent V\'eron}
\vspace{3mm}\\
{\it Laboratoire de Math\'ematiques et Physique Th\'eorique, CNRS UMR 6083},
\\
{\it Universit\'e Fran\c{c}ois Rabelais},{\it Tours 37200, France}}

\date{}
\maketitle
\begin{center} {\bf\small Abstract}\end{center}{\small  If $h$ is a nondecreasing real
valued function and
$0\leq q\leq 2$, we analyse the boundary behaviour of the gradient of any solution $u$ of 
$-\Delta u+h(u)+\abs {\nabla u}^q=f$ in a smooth N-dimensional domain $\Gw$ with the condition that $u$  tends to
infinity when $x$ tends to $\prt\Gw$. We give precise expressions of the blow-up which, 
in particular, point out the fact that the phenomenon occurs essentially in the normal direction to
$\prt\Gw$.  Motivated by the blow--up argument in our proof, we also give in Appendix  a symmetry   result  for  some related problems in the half space.}

\section {Introduction}

Let $\Gw$ be a $C^2$ domain in $\rn$ ($N\geq 2$), $h$ a continuous nondecreasing function and $q$ a
nonnegative real
number. The
aim of this
work is to study the behaviour of solutions of nonlinear equations of the following type
\be\label{equa} -\Delta u+h(u) +|\nabla u|^q = f\qquad\hbox{ in $\Omega\subseteq \rn$,}
\ee satisfying  a boundary blow--up condition
\be\label{bc}
\lim\limits_{\dom \to 0} u(x)=+\infty\,
\ee 
where $d_\Gw(x)=\dist (x,\prt\Gw)$. The interest for solutions of \rife{equa} satisfying
such singular boundary conditions
arises from stochastic control problems with state constraints, as explained in \cite{LL}, where $h(u)=\la\, u$. In
that situation, $u$ represents the value function of the optimal control problem  and $ -q\nabla u\,|\nabla
u|^{q-2}$ acts as the optimal  (feedback) control which forces 
the process to stay in  $\Omega$.  
 
From  a purely PDE's point of view, the existence of such solutions depends on   the
possibility  of finding universal interior estimates for \rife{equa}, independently on the
behaviour of $u$ at the boundary. In the case $q=0$ these estimates hold provided the
well--known  Keller--Osserman  condition (\cite{K}, \cite{O}) is satisfied, i.e.
\be\label{OK}
 \int^{+\infty} \frac{ds}{\sqrt{\int_0^sh(t)dt}}<\infty\,.
\ee A  large number of papers has investigated properties of such singular solutions (also
called {\it large}, or {\it explosive} solutions) 
when the lower order terms only depend on $u$ (see \cite{BaM}, \cite{BaM2}, \cite{D-L},
\cite{LN}, \cite{MV}, \cite{MV1}, \cite{Ve}).  In presence of gradient dependent terms  
as in \rife{equa},  large solutions in smooth domains have been studied in \cite{BG},
\cite{G},  \cite{GNR}, \cite{LL}, \cite{Ploc}; roughly speaking,  such solutions exist
if   $h$ satisfies \rife{OK} or if $1<q\leq 2$ and $h$ is unbounded at infinity. Indeed,
in equation \rife{equa} both lower order terms may lead to the construction of large
solutions, so that existence of solutions to problem \rife{equa}--\rife{bc} can be proved
even if $h$ is sublinear, provided $q>1$.

In this paper we consider problem  \rife{equa}--\rife{bc}, mainly referring to the model
examples $h(s)=e^{as}$, $a>0$, and $h(s)=s^\beta$, $\beta> 0$, and we study the asymptotic behaviour of 
$\nabla u$ at the boundary. 
It turns out, as a quite general rule, that $\nabla u$ blows up, in its first
approximation, in the normal direction: in  the model examples, our results read as
follows. We denote by 
$d_{ \Omega}(x)$ the distance of a  point $x$ to $
\partial \Omega$, and  by $\nu$ the outward unit normal  vector at $
\partial \Omega$.

 \begin{theorem}\label{tex}
 Let $\Omega$ be a $C^2$ domain in $\rn$,  $\nu$ be the normal outward unit vector  to $
\partial
\Omega$, and assume  $f\in \elle\infty$. \smallskip 

\noindent{\bf A-} Let  $a> 0$, and $u$ be a solution of 
 $$
 \cases{-\Delta u+ \, e^{au} +|\nabla u|^q = f& in $\Omega$, \cr \lim\limits_{\dom\to
0}u(x) =+\infty\,.\cr} $$
Then there holds:
 \begin{itemize}
\item[(1)] If $q=2$ and $a\leq 2$,    then 
$$
 \lim\limits_{\dom\to 0} \dom \nabla u(x)=\nu. $$
 \item[(2)] if $0\leq q< 2$, or if $q=2$ and $a>2$,  then  
$$
 \lim\limits_{\dom\to 0} \dom \nabla u(x)=\frac 2a\nu. $$
\end{itemize}\smallskip

\noindent {\bf B-} Let   $\beta >0$ and $u$ be a    solution of 
$$
 \cases{-\Delta u+ |u|^{\beta-1}u +|\nabla u|^q = f& in $\Omega$, \cr \lim\limits_{\dom\to
0}u(x) =+\infty\,.\cr} $$
Then there holds:
 \begin{itemize}
\item[(3)] If  $q\geq  \frac{2\beta}{1+\beta}$,  then   
$$
\lim\limits_{\dom\to 0}\dom^{\frac1{q-1}}\, \nabla u(x)=b\,\nu,$$
\end{itemize}
in which formula
$b= (q-1)^{-\frac1{q-1}} $ 
if $q>  \frac{2\beta}{1+\beta}$, and
$b=\left(\frac1a\right)^{\frac{2-q}{2(q-1)}}\left(\frac{2-q}{q-1}\right)^{\frac1{q-1}}$
 if
$q=\frac{2\beta}{1+\beta}$, where $a$ is the solution of $a-a^{\frac q2}= 2-q $.
\begin{itemize}
\item[(4)] If  $q< \frac{2\beta}{1+\beta}$,  then 
$$
 \lim\limits_{\dom\to 0} \dom^{\frac{1+\beta}{\beta-1}} \nabla u(x)=b\, \nu,$$
 \end{itemize}
where 
$b=\frac2{\beta-1}\left[\frac{2(\beta+1))}{(\beta-1)^2)}\right]^{1/(\beta-1)}$.
 \end{theorem}
 
The previous result generalizes those obtained in \cite{BE} and \cite{BaM2} for large
solutions
 of semilinear problems, in case the lower order terms do not depend on $\nabla u$;
indeed, our proof   follows a  similar approach based on  a blow--up argument near the
boundary and requires  some symmetry results on the blown--up functions,  which are
solutions  of   a similar problem in the half space.  Even in the case $q=0$, our result
extends those previous ones by considering  a slightly larger class of nonlinearities $h(s)$. The
conclusions  of Theorem \ref{tex} will follow as a particular case of the results which we
prove in Section 2. 
Moreover, in  a third section    we will also provide  a  simple uniqueness result for solutions of
\rife{equa}--\rife{bc}
 which is meant to be applied in case  $h$ is concave, or the sum of a  concave and a
convex function.
 In fact,  previous uniqueness results  
   seem to have been proved only if $h$ has a  convex type behaviour.
   
Finally, motivated by our blow--up  argument in case $h(s)$ has a power growth at infinity, we prove in Appendix   some symmetry and uniqueness results for  nonnegative solutions of     
the problem in the half space 
$$
\cases{-\Delta u+\al\,  u^{p}+ |\nabla u|^q =0&in $\rnp\,:\,=\{\xi=(\xi_1,\xi')\in \rn\,:\,\xi_1>0\}$,\cr
\noalign{\medskip}
u(0,\xi')=M & \cr}
$$
 where $\alpha\geq 0$, $p>0$ and $M$ is a  nonnegative constant or possibly $M=+\infty$.
We give  a simple proof, based mainly on comparison with radial or one--dimensional solutions, that any nonnegative solution $u$
is one--dimensional, and uniqueness follows if $\alpha>0$.
 \section{Asymptotic behaviour of derivatives}
 
 In this section we let $\Omega\subset \rn$ be a  bounded $C^2$ domain. We denote by 
$d_{ \Omega}(x)=\dist(x,
\partial \Omega)$, and  by $\nu(x)$ the outward unit normal  vector at any point $x\in 
\partial \Omega$, or  simply $\nu$ when meant as a vector field defined on $
\partial \Omega$. In the sequel,   $\tau$ is any unitary tangent vector field defined on
$\partial \Omega$ as well, i.e. $\tau\cdot \nu=0$.

We start by considering the equation 
 \be\label{q2}
 \cases{-\Delta u+h(u) +|\nabla u|^2 = f& in $\Omega$, \cr \lim\limits_{\dom\to 0}u(x)
=+\infty\,,\cr}
 \ee where $h$ is an increasing function such that $\lim\limits_{s\to
+\infty}h(s)=+\infty$, and $f\in \elle\infty$. 

It is proved in \cite{Ploc}  that problem
\rife{q2} admits a solution, and moreover any solution satisfies the estimate
\be\label{por}
 \hbox{$u(x)- F(\dom)$  is bounded near $
\partial \Omega$, where $F^{-1}(s)=\dys \int^{+\infty}_s \frac{e^{-t}}{[ \int_0^t
h(\xi)e^{-2\xi}d\xi]^{\frac12}}dt$.}
\ee Note  that the function $F$ has at most  a logarithmic blow--up rate. Moreover, if the
following limit exists $$
\lim\limits_{\xi \to +\infty}\left(1+ \frac12 \frac{h(\xi)e^{-2\xi}}{\int_0^\xi
h(t)e^{-2t}dt}\right)^{-1} $$
 one has,  using twice  L'Hopital's rule and  since both $F^{-1}(\xi)$ and $(F^{-1})'(\xi)$
tend to zero as $\xi$ goes to infinity,
\be\label{dlog}
\eqalign{ &\qquad \qquad 
\lim\limits_{s\to 0} \frac{F(s)}{|\log s|} = 
-\lim\limits_{s\to 0} s\, F'(s)= \cr & =
-\lim\limits_{\xi \to +\infty}  \frac{F^{-1}(\xi)}{(F^{-1})'(\xi)}=-\lim\limits_{\xi \to
+\infty} \frac{(F^{-1})'(\xi)}{(F^{-1})''(\xi)}=
\lim\limits_{\xi \to +\infty} \left(1+ \frac12 \frac{h(\xi)e^{-2\xi}}{\int_0^\xi
h(t)e^{-2t}dt}\right)^{-1}\,.\cr}
\ee
 Similarly one has
\be\label{+log}
\eqalign{ &\qquad \qquad 
\lim\limits_{s\to 0} \left(F(s)+\log s\right)=\lim\limits_{\xi\to +\infty} \log(e^\xi F^{-1}(\xi))=\cr
&= \log\left(-\lim\limits_{\xi\to+\infty}\frac{(F^{-1})'(\xi)}{e^{-\xi}}\right)=-\frac12 \log
\left(\lim\limits_{\xi \to +\infty}\int_0^\xi h(t)e^{-2t}dt\right)\,.\cr}
\ee 
In particular we deduce that
\be\label{subc}
\hbox{$u(x)+\log(\dom)$ is bounded near $
\partial \Omega$ if and only if $\displaystyle{\int_0^{+\infty} h(t)e^{-2t}dt<\infty}$,}
\ee and that
\be\label{supc}
\hbox{if $\dys\lim\limits_{s\to +\infty} \frac{h(s)e^{-2s}}{\int_0^s h(t)e^{-2t}dt}=\la\geq
0$, then $\dys\frac{u(x)}{|\log(\dom)|}\to \frac2 {\la+2}$ as $\dom \to 0$.}
\ee

In view of these remarks,  we will consider three types of situations in our analysis,
which are mutually excluding:

\begin{itemize}
\item[(h1)] $\dys \int^{+\infty} h(t)e^{-2t}dt <\infty$ and $\lim\limits_{s\to +\infty}
h(s)e^{-2s}=0$.
\item[(h2)] $\dys\int^{+\infty} h(t)e^{-2t}dt =\infty$, $\dys\lim\limits_{s\to +\infty}
\frac{h(s)e^{-2s}}{\int_0^s h(t)e^{-2t}dt}=0$, and $\dys\frac{h(s+c)}{h(s)}$ is bounded for
large $s$, and any $c\in \r$.
\item[(h3)] $\dys\lim\limits_{s\to +\infty} \frac{h(s)e^{-2s}}{\int_0^s h(t)e^{-2t}dt}=\la>0$,
and, for any $t\in \r$, $\exists\dys\lim\limits_{s\to +\infty} 
\frac{h(s+t)}{h(s)}=e^{(\la+2)t}$.
\end{itemize}
\bigskip

\begin{remark}\label{tre}\rm Assumption (h1) corresponds to  a subcritical case, where the
  blow--up rate of $u$ only depends on the first order term, whereas  (h2) represents the
critical case (e.g.  $h(s)=e^{2s}$) in which both terms give  a contribution and a
superposition effect may be observed; in fact, due to \rife{subc}--\rife{supc}, in both
cases we have $\dys\frac{u(x)}{|\log(\dom)|}\to 1$, but while under (h1) we have that
$u(x)+\log(\dom)$ is bounded near $\partial \Omega$, (h2) implies that $u(x)+\log(\dom)\to
-\infty$ at the boundary.

As far as (h3) is concerned, it covers exponential--type growths, including the model
$h(s)=e^{(2+\la)s}s^\beta$ for any $\beta\geq 0$. Let us remark that  assuming the
existence,  for any $t\in \r$, of $\dys\lim\limits_{s\to +\infty}  \frac{h(s+t)}{h(s)}$  {\sl
automatically} implies that the function 
$\omega(t)\,:\,=\dys\lim\limits_{s\to +\infty}  \frac{h(s+t)}{h(s)}$ is an exponential.
Indeed, since $h$ is increasing, the same is true for  $\omega $. Since   
$\omega(t+t')=
\omega(t)\omega(t')$ for every  $t$, $t'\in \r$, the continuity of $\omega$ at  a point $t_0$ implies that 
$\omega $ is continuous on $\r$, and then (using also  
$\omega(0)=1$) $\omega(t)=e^{a\,t}$ for some $a\in \r$.  Moreover, since $\omega$ is continuous 
the above convergence is locally uniform for $t$ in $\r$.
Eventually, if 
\be\label{nuo}
\la= \lim\limits_{s\to +\infty} \frac{h(s)e^{-2s}}{\int_0^s h(t)e^{-2t}dt} \,, 
\ee
we have 
$$
\frac{\int_0^{s+t}h(\xi)e^{-2\xi}d\xi}{e^{-2s}h(s)}= 
\frac{\int_0^{s}h(\xi)e^{-2\xi}d\xi}{e^{-2s}h(s)}+\int_0^t \frac{
h(s+\xi)}{h(s)}e^{-2\xi}d\xi\to \frac1\la +\int_0^t
e^{(a-2)\xi}d\xi
$$
as $s\to\infty$. But L'Hopital's  rule also implies 
$$
\lim\limits_{s\to
+\infty}\dys\frac{\int_0^{s+t}h(\xi)e^{-2\xi}d\xi}{\int_0^{s}h(\xi)e^{-2\xi}d\xi}= e^{(a-2)t},
$$
so that we deduce, using also 
\rife{nuo},
 $$
\frac1\la+\int_0^t e^{(a-2)\xi}d\xi=
\lim\limits_{s\to
+\infty}\frac{\int_0^{s+t}h(\xi)e^{-2\xi}d\xi}{e^{-2s}h(s)}=\frac{e^{(a-2)t}}{\la}, $$
hence $a\neq 2$, and $a=\la+2$.
\end{remark}
\qed
\vskip2em

 \begin{theorem}\label{tq2}
 Let $u$ be a solution of \rife{q2}. Then we have:
 \begin{itemize}
\item[(1)] If $(h1)$ or (h2) hold true,
\be\label{ldir}
\lim\limits_{\de \to 0}\de\,\frac {
\partial u}{
\partial \nu(x)}(x-\de\nu(x))=1\,,\qquad \lim\limits_{\de \to 0}\de\, \frac{
\partial u}{
\partial \tau(x)}(x-\de\nu(x))=0
\ee
 \end{itemize}
holds uniformly for $x\in
\partial\Omega$, and then 
\be\label{lglo}
 \lim\limits_{\dom\to 0} \dom \nabla u(x)=\nu.
\ee
 \begin{itemize}
 \item[(2)] If (h3) holds true, 
\be\label{ldir2}
\lim\limits_{\de \to 0}\de\, \frac{
\partial u}{
\partial \nu(x)}(x-\de\nu(x))=\frac2{\la+2}\,,\qquad \lim\limits_{\de \to 0}\de\, \frac{
\partial u}{
\partial \tau(x)}(x-\de\nu(x))=0
\ee 
\end{itemize}
holds uniformly for $x\in \partial\Omega$, and then 
\be\label{lglo2}
 \lim\limits_{\dom \to 0} \dom \nabla u(x)=\frac 2{\la+2} \nu.
\ee
 \end{theorem}
 \Proof Thanks to \rife{por}, we can fix $d_0$ and $C_0$ such that 
\be\label{porr}
 \eqalign{ & |u(x)-F(\dom)|\leq C_0\quad \hbox{for any $x\in \Omega$: $\dom \leq d_0$,}
\cr &\quad \hbox{ with $F^{-1}(s)= \dys\int^{+\infty}_s \frac{e^{-t}}{[\int_0^t
h(\xi)e^{-2\xi}d\xi]^{\frac12}}dt$.}\cr}
 \ee
 We use  a similar blow--up framework as in \cite{BE},  \cite{BaM2}.
 Let $x\in 
\partial \Omega$ and consider a new system of coordinates 
$(\eta_1,  \ldots, \eta_N)$ centered at $x$ and such that the positive $\eta_1$-axis is 
the direction $-\nu(x)$, where $\nu(x)$ is the outward normal vector at $x$; thus $x=O$ is
the origin and $\eta_1$ is the direction of the inner normal vector at $x$. In the
$\eta$--space, let us set $P_0=(d_0,0,\ldots,0)$  and define
 $$
 D_\de= B(O,\de^{1-\sigma})\cap B(P_0, d_0)\,,\qquad \hbox{with $0<\sigma<\frac 12$.}
 $$
Note that we can assume that $\Gw$ satisfies the interior sphere condition with radius $d_{0}$ so that $D_\de\subset \Omega$, and since the operator is invariant under translations
and 
rotations we obtain the same equation for $u$ in the new variable $\eta$. Define
$\xi=\frac{\eta}\de$ and the function
 $$
 v_\de(\xi)= u(\eta)- F(\de)= u(\de\xi)-F(\de)\,,
 $$
 where $F$ is defined in \rife{porr}. Then $v_\de(\xi)$ satisfies the equation
 $$
 -\Delta v_\de +h(u(\de\xi)) \de^2+|\nabla v_\de|^2=\de^2 f(\de\xi)\qquad \xi \in
\frac1\de
 D_\de.
 $$
 It is readily seen that since  $0<\sigma<\frac12$,  if $\eta\in \partial B(P_0, d_0)\cap
\partial D_\de$, then $\frac{\eta_1}\de\to 0$ and $\frac{|\eta'|}\de\to +\infty$ as
$\gd\to 0$;
moreover
 since $|\eta|<\de^{1-\sigma}$, we conclude that the domain $ \frac1\de D_\de$
converges to the half space
 $\rnp\,:\,=\{\xi\in \rn\,:\, \xi_1>0\}$.
 
 Let us study now the limit  of $v_\de$.  First of all, observe that  
 since $F^{-1}$ is a  decreasing and convex function (as easily checked), then its inverse
 function $F$ is also convex. We have then, for any $\la<1$, 
 $$
0\leq  F(\la s)-F(s) \leq -F'(\la s)\,\la s\,\,\frac{1-\la}{\la}\,,
 $$
 and since (see also \rife{dlog}) $0 < -F'(\xi)\xi< C$ for any $\xi\in \r^+$, we deduce
that
  $F$   enjoys the property
 \be\label{prof}
\exists C>0\,:\quad  F(\la s)-F(s) \leq C\frac{1-\la}{\la}\qquad \forall \la<1\,, \quad
\forall s>0\,.
 \ee
 Since $\partial\Omega$ is $C^2$, we have that  for $\eta \in D_\de$
\be\label{Tay}
 d_{\Omega}(\eta)= \eta_1+ O(|\eta|^2)= \eta_1+ O(\de^{2-2\sigma})\,. 
 \ee
 Hence 
 from  \rife{porr}--\rife{prof} we deduce that
 \be\label{alw}
 |u(\de\xi)- F(\de\,\xi_1+\de^{2-2\sigma})| \leq C_1\quad \hbox{for any $\xi\in \frac1\de
D_\de$,}
 \ee so that 
 \be\label{estv} | v_\de(\xi)|\leq C_1 +|F(\de\,(\xi_1+\de^{1-2\sigma}))- F(\de)| \quad
\hbox{for any $\xi\in \frac1\de D_\de$.}
\ee
 In particular, due to \rife{prof}, \rife{estv} implies that
 $$
 |v_\de(\xi)|\leq C_1 + C_2 \,\max\{\xi_1\,,\,\frac1{\xi_1}\}\,, 
 $$
 hence $v_\de $ is   locally uniformly bounded. 
 
Assume that (h1) holds true: then (see \rife{+log})  
$F(\de)+\log(\de)$ is bounded for small $\de$, so that \rife{estv} implies that $$
v_\de(\xi)\geq F(\de(\xi_1+\de^{1-2\sigma}))-F(\de)- C_1 \geq -\log(\xi_1+\de^{1-2\sigma})
-C_2\,, $$
for $\xi \in \frac1\de D_\de$; in particular in the limit (as $\de\to 0$)  we deduce (recall that
$\sigma<\frac12$)
\be\label{dab} v(\xi)\geq  -\log\xi_1 -C_2
\ee so that $\lim\limits_{\xi_1\to 0^+} v(\xi)=+\infty$.
Noticing that 
 $$
\de^2 h(u(\de\xi))= h(v_\de+F(\de))e^{-2(v_\de+F(\de))}\, e^{2(v_\de+ F(\de)+\log\de)}\leq
C h(v_\de+F(\de))e^{-2(v_\de+F(\de))}\, e^{2v_\de}\,, $$
and using that $v_\de$ is locally bounded and $h(s)e^{-2s}\to 0$ as $s\to +\infty$, we deduce
\be\label{lot}
\de^2 h(u(\de\xi))\to 0 \quad \hbox{in $L^\infty_{loc}(\rnp)$}\,.
\ee Furthermore, standard  elliptic estimates   for second derivatives imply that $|\nabla
v_\de|$ is also 
 locally uniformly bounded, and, in the end, that   $v_\de$ is locally relatively compact
in the  $C^1_{loc}$--topology.
 Let $v$ be the limit of some subsequence $v_{\de_k}$, as $\de_k\to 0$. Therefore $v$ is a
solution of
\be\label{eqv}
\cases{ -\Delta v+|\nabla v|^2=0 & in $\rnp$, \cr \lim\limits_{\xi_1\to 0^+}v(\xi)
=+\infty\,.\cr}
\ee The function $w=e^{-v}$ is positive and harmonic in $\rnp$; it satisfies $w\leq C\xi_1$, from \rife{dab}, hence $w=0$ on $\{\xi_1=0\}$. We deduce (for instance
using Kelvin transform, or symmetry results) that there exists $\la\in \r_+$ such that
$w=\la\,\xi_1$, hence $v=-\log\xi_1-\log \la$. In particular, we obtain, locally uniformly
in $\rnp$: $$
\frac{
\partial v_{\de_k}}{
\partial \xi_1}\to -\frac1{\xi_1}\,,\qquad \frac{
\partial v_{\de_k}}{\partial
\xi_j}\to 0\quad \forall j=2,\ldots,N, $$
for any convergent subsequence $v_{\de_k}$. Note that while the limit function $v$ is
determined up to the  constant $-\log \la$, its gradient is uniquely  determined. This
implies  that the whole sequence of derivatives $\frac{\partial v_{\de}}{\partial
\xi_i}$ will be converging to this limit. We have proved then that it holds: $$
\de\frac{
\partial u(\de\xi)}{
\partial \xi_1}\to -\frac1{\xi_1}\,,\qquad \de\frac{
\partial u(\de\xi)}{\partial
\xi_j}\to 0\quad \forall j=2,\ldots,N. $$
Recalling that $\xi_1$ is the direction of the inner normal vector and that the point
$\eta=(\de, 0, \ldots,0)$ coincides with  $x-\de\nu(x)$, we fix $\xi_1=1$ and obtain
\rife{ldir}.
 
Let us now assume (h2). In this case $F(\de)+\log(\de)$ is unbounded,  but we still have
(see \rife{dlog}) $$
F'(\de)\de\to -1\qquad\hbox{as $\de\to0$.} $$
In particular, for any $\ga<1$ there exists  an interval $(0, s_\ga)$  such that the
function $F(s)+\ga\log s$ is decreasing in $(0,s_\ga)$; therefore, for $\xi_1<1$  and
$\de$ small enough, we have 
$$
F(\de(\xi_1+\de^{1-2\sigma}))-F(\de)\geq -\ga\log(\xi_1+\de^{1-2\sigma})\,. $$
Together with \rife{estv} we deduce that $$
v_\de(\xi)\geq F(\de(\xi_1+\de^{1-2\sigma}))-F(\de)- C_1
\geq-\ga\log(\xi_1+\de^{1-2\sigma})-C_1 $$
hence, for any possible limit function $v$, we deduce that $v\geq -\ga\log\xi_1-C_1$ for
$\xi_1$ near zero. This implies in particular that $v$ blows--up uniformly on $\{\xi_1=0\}$.
Writing again
\be\label{hd}
\de^2 h(u(\de\xi))= \frac{h(v_\de+F(\de))}{h(F(\de))}\,
\frac{h(F(\de))e^{-2F(\de))}}{\int_0^{F(\de)}h(s)e^{-2s}ds}\, e^{2\log (\de
e^{F(\de)}[\int_0^{F(\de)}h(s)e^{-2s}ds]^{\frac12})}\,,
\ee 
and using (h2) and (see \rife{dlog}) 
$$
\lim\limits_{t\to+\infty}
F^{-1}(t)e^t[\int_0^{t}h(s)e^{-2s}ds]^{\frac12}=\lim\limits_{t\to+\infty}
-\frac{F^{-1}(t)}{(F^{-1})'(t)}=1, $$ 
we conclude that \rife{lot} still holds true. Then, passing to the limit in $\de$, any
limit function $v$ will satisfy \rife{eqv}. Again, we have that $w=e^{-v}$ is  harmonic in
$\rnp$ and $w\leq C \xi_1^\ga$ in a  neighborhood of $\{\xi_1=0\}$, so that $w=0$ on
$\partial \rnp$. We conclude as above that $w=\la \xi_1$ for some $\la \in \r_+$, and then
$v=-\log\xi_1-\log \la$. As before, the convergence of $\nabla v_\de$ to $\nabla v$  then 
implies \rife{ldir} and \rife{lglo}.

Finally, let us assume (h3), and let again $v$ be such that (a subsequence of) $v_\de$
converges to $v$ locally uniformly. Due to the monotonicity of $h$, we have (see Remark
\ref{tre}): $$
\lim\limits_{s\to +\infty}\frac{h(s+t)}{h(s)}=e^{(\la+2)t}\quad \hbox{locally uniformly in
$t$} $$
so that  
$$
\lim\limits_{\de\to0} \frac{h(v_\de+F(\de))}{h(F(\de))}=e^{(\la+2)v}\quad\hbox{in
$L^\infty_{loc}({\rnp})$}. $$
Since under (h3) we also have  (see \rife{dlog})
\be\label{che}
\lim\limits_{t\to+\infty}
F^{-1}(t)e^t[\int_0^{t}h(s)e^{-2s}ds]^{\frac12}=\lim\limits_{t\to+\infty}
-\frac{F^{-1}(t)}{(F^{-1})'(t)}= \lim\limits_{s\to 0} - F'(s)s=\frac2{\la+2},
\ee
 then  \rife{hd} now implies
\be\label{hlim}
\lim\limits_{\de\to 0}
\de^2 h(u(\de\xi))= e^{(\la+2)v}\,\la\, e^{2\log (\frac2{\la+2})}= c_\la e^{(\la+2)v}
\ee where $c_\la=\frac{4\la}{(\la+2)^2}$. Moreover we also deduce   from \rife{che}  that
 there exist an interval $(0,\sigma_0)$ and constants $\ga_0<\frac 2{\la +2}$
 and $\ga_1>\frac 2{\la +2}$ such that $F(t)+\ga_0\log t$ is decreasing and
$F(t)+\ga_1\log t$ is increasing  in $(0,\sigma_0)$. In particular we have 
$$
F(\de(\xi_1+\de^{1-2\sigma}))-F(\de)\geq -\ga_0\log( \xi_1+\de^{1-2\sigma}) \qquad
\hbox{if $\xi_1\leq 1-\de^{1-2\sigma}$,} 
$$
and $$
F(\de(\xi_1+\de^{1-2\sigma}))-F(\de)\geq -\ga_1\log( \xi_1 +\de^{1-2\sigma})\qquad
\hbox{if $1<\xi_1< \frac{\sigma_0}\de-\de^{1-2\sigma}$,} $$
which together with \rife{estv} imply
\be\label{zer} v_\de(\xi)\geq -\ga_0\log( \xi_1 +\de^{1-2\sigma})- c_0\qquad \hbox{if
$\xi_1\leq 1-\de^{1-2\sigma}$,}
\ee and
\be\label{inf} v_\de(\xi)\geq -\ga_1\log( \xi_1 +\de^{1-2\sigma}) -c_1\qquad \hbox{if
$1<\xi_1< \frac{\sigma_0}\de-\de^{1-2\sigma}$.}
\ee From \rife{hlim} and \rife{zer}--\rife{inf} we deduce, passing to the limit in $\de$, 
that $v$ satisfies
\be\label{eqv2}
\cases{ -\Delta v+c_\la\,e^{(\la+2)v}+|\nabla v|^2=0 & in $\rnp$, \cr
\lim\limits_{\xi_1\to 0^+}v(\xi) =+\infty\,,&\cr}
\ee and the further estimate
\be\label{add} v(\xi)\geq -\ga_1\log \xi_1 -c_1\qquad \hbox{if $1<\xi_1$.}
\ee 
We proved in \cite{PoVe} (Corollary 2.6) that any solution of \rife{eqv2} only depends
on the $\xi_1$ variable, moreover condition  \rife{add} implies that we have exactly $$
v=\frac2{\la+2}\log(\frac1{\xi_1})+
\frac1{\la+2}\log(\frac{2\la}{c_\la(\la+2)^2})=
\frac2{\la+2}\log(\frac1{\xi_1})-
\frac{\log 2}{\la+2}\,. $$
We obtain that $$
\frac{
\partial v_\de}{
\partial \xi_1}\to -\frac2{(\la+2)\xi_1}\,,\qquad \frac{
\partial v_\de}{\partial
\xi_j}\to 0\quad \forall j=2,\ldots,N, $$
which, as before, gives \rife{ldir2} and \rife{lglo2}.
\qed

\begin{remark} \rm The same proof applies if one  only requires on the right hand side  that
$\lim\limits_{\dom \to 0}d^2_\Gw(x) f(x)=0$, which implies
 that $\lim\limits_{\de\to 0} \de^2 f(\de\xi)= 0$ locally uniformly for  $\xi\in\rnp$.
\end{remark}

\begin{remark} \rm Under assumption (h3),  the previous proof gives that the rescaled sequence
$v_\de$ converges towards $v=\frac2{\la+2}\log(\frac1{\xi_1})-
\frac{\log 2}{\la+2}$. Setting $\xi_1=1$ we deduce that 
$$
u(x)-F(\dom)\to -
\frac{\log 2}{\la+2} $$
which improves estimate \rife{por}. As a consequence, this also implies that
$u_1(x)-u_2(x)\to 0$ for any two large solutions $u_1$, $u_2$, hence in this case
uniqueness of solutions of \rife{q2} follows immediately by the maximum principle.
\end{remark}

\vskip1em We consider now the problem
\be\label{q}
 \cases{-\Delta u+h(u) +|\nabla u|^q = f& in $\Omega$, \cr \lim\limits_{\dom\to 0}u(x)
=+\infty\,,&\cr}
 \ee with $0\leq q<2$. In this case if $h$ has an exponential growth at infinity, the gradient
term does not affect the behaviour of solutions near the boundary, so that the asymptotic
behaviour of this problem turns out to be  the same as for the semilinear equation with 
$q=0$. In order to  adapt  the above proof we will need the following uniqueness result
for solutions in the half space.

\begin{lemma}\label{sym} Let $a>0$ and $v$ be a solution of 
$$
\cases{-\Delta v+ \,e^{av}=0 & in $\rnp$, \cr \lim\limits_{\xi_1\to 0^+}v(\xi) =+\infty&
locally uniformly with respect to $\xi'\in {\bf R}^{N-1}$\,.\cr} $$
Assume that $v$ satisfies the following assumption:
\be\label{asgro}
\exists \gamma\,,\,m\,,\,S_0>0\,:\quad v(\xi)\geq -\gamma\log S-m\quad \forall
\xi\in\rn\,:\,\xi_1\leq S\,,\quad \forall S> S_0\,.
\ee Then $v= -\frac2a\log\xi_1+\frac1a\log\frac2a$.   
\end{lemma}
\Proof We can assume $a=1$, up to replacing $v$ with $\frac 1a v- \frac 1a \log a$.  We follow the
approach  used in \cite{PoVe}  (see Proposition 4.1); for any $R>0$, $S>S_0$, define
$\omega_R$  as the solution of the  problem $$
\cases{ -\Delta \omega_R+ \,e^{\omega_R}=0 & in $B_R(0)$, \cr
\lim\limits_{\rho\uparrow R}\omega_R(\rho) =+\infty\,,&\cr} $$
and define ${\underline \omega}_{R,S}$ as the solution of the problem 
$$
\cases{ -\Delta {\underline \omega}_{R,S}+ \,e^{{\underline \omega}_{R,S}}=0  & in $B_{R+S}(0)\setminus B_R(0)$, \cr
\lim\limits_{\rho\downarrow R}{\underline \omega}_{R,S}(\rho) =+\infty\,,\quad &$
{\underline \omega}_{R,S}(R+S)=-\ga\log S-m$.\cr} $$
Now fix $\xi'\in {\bf R}^{N-1}$, and consider  the points $\xi_R=(R, \xi')$, 
$\eta_R=(-R, \xi')$  and the functions $\omega_R(\cdot-\xi_R)$ and
${\underline \omega}_{R,S}(\cdot-\eta_R)$. By comparison, 
and using \rife{asgro}, we have  
\be\label{comp} v\leq \omega_R(\cdot-\xi_R) \quad \hbox{in $B_R(\xi_R)$,}\qquad v\geq
{\underline \omega}_{R,S}(\cdot-\eta_R)\quad \hbox{in 
$B_{R+S}(\eta_R)\cap \rnp$.}
\ee It is readily seen that the sequence $\{\omega_R(\cdot-\xi_R)\}$ is decreasing 
and converges, as $R\to +\infty$, to  a  function $\omega_{\infty}$ which only depends
on the $\xi_1$--variable and is  the maximal solution of 
\be\label{1d} -z''+e^{z}=0\,,\quad \lim\limits_{t\to 0^+}z(t)=+\infty\,.
\ee In particular, from a  straightforward computation of solutions of \rife{1d}, we obtain
$\omega_{\infty}(\xi_1)=-2\log\xi_1+\log2$.

Let $S>S_0$; without loss of generality we can replace the constants $\ga$ and $m$ in
\rife{asgro} with possibly larger 
values. In particular, we can assume that  $\ga>2$ and $e^{-m}<2S_0^{\ga-2}$: let then
$w(\rho)= -2\log(\rho-R)-(\ga -2)\log S-m$,  computing we have, for $\rho\in (R,R+S)$: $$
\eqalign{ 
-\Delta w+e^w&=
\frac{2(N-1)(\rho-R)S^{\ga-2}-(2S^{\ga-2}-e^{-m})\rho}{(\rho-R)^2S^{\ga-2}\rho}\cr & \leq
\frac{2(N-1) S^{\ga-1}-(2S^{\ga-2}-e^{-m})R}{(\rho-R)^2S^{\ga-2}\rho}\,,\cr} $$ 
so that there exists a value $R_0(S)$ such that 
$$
-\Delta w+e^w\leq 0 \qquad \hbox{in $B_{R+S}(0)\setminus B_R(0)$ for any $R\geq R_0(S)$.}
$$
Since $w(R+S)=-\ga\log S-m$ we deduce that 
$$
{\underline \omega}_{R,S}\geq w \geq -\ga\log S-m \qquad \hbox{for any $R\geq R_0(S)$.} $$
In particular, for any $R>R' >R_0(S)$, comparing ${\underline \omega}_{R,S}(\cdot-\eta_R)$ and
$\underline\omega_{R',S}(\cdot-\eta_{R'})$ (on their common domain $B_{R'+S}(\eta_{R'})\setminus
B_R(\eta_R)$) we deduce that 
$$
{\underline \omega}_{R,S}(\cdot-\eta_R)\geq \omega_{R',S}(\cdot-\eta_{R'}) $$
hence for any fixed $S$ the sequence $\{{\underline \omega}_{R,S}(\cdot-\eta_R)\}_{R}$ is definitively
increasing and converges to 
a  function ${\underline \omega}_{S}$ which only depends on the $\xi_1$--variable and solves
\be\label{1ds} -{\underline \omega}_{S}''+e^{{\underline \omega}_{S}}=0\,,\quad \lim\limits_{t\to
0^+}{\underline \omega}_{S}(t)=+\infty\,,\quad {\underline \omega}_{S}(S)= -\ga \log S-m\,.
\ee Thus from \rife{comp}, passing to the limit in $R$, we derive
\be\label{quas}
{\underline \omega}_{S}(\xi_1)\leq v(\xi)\leq -2\log\xi_1+ \log 2\qquad \forall
\xi\in\rnp\,:\,\xi_1\leq S\,,\quad \forall S>S_0\,.
\ee 
Next, letting $e^{-m}\leq 2$, we observe that the function $z$ defined by $z(t)=-2\log
t-(\ga-2)\log(t+1)-m$ satisfies $$
-z''+e^z=- \frac2{t^2}-\frac{\ga-2}{(t+1)^2}+\frac{e^{-m}}{t^2(t+1)^{\ga-2}}\leq
\frac{-2(t+1)^{\ga-2}+e^{-m}}{t^2(t+1)^{\ga-2}}\leq 0\,, $$
and since $z(S)<-\ga\log S-m$ we have that it is a subsolution for the problem \rife{1ds},
hence 
\be\label{logga} -2\log t-(\ga-2)\log(t+1)-m\leq{\underline \omega}_{S}(t)\leq -2\log t+ \log 2\,.
\ee The sequence $\{{\underline \omega}_{S}(t)\}_{S\geq S_0}$ is then locally bounded and, up to
subsequences,  converges (locally in the $C^2$--topology) to  a solution
${\underline \omega}_{\infty}$ of \rife{1d}; but  estimate \rife{logga}  implies (due to  the
classification of all solutions of \rife{1d}, see e.g. \cite{PoVe})  that the only
possible limit  is ${\underline \omega}_{\infty} = -2\log t+\log 2$. Letting $S$ go to infinity, we
conclude from \rife{quas} that $v=-2\log\xi_1+ \log2$.
\qed

We are ready now to deal with the case that $q<2$ and $h$ has an  exponential scaling at
infinity. Our next result extends the one in \cite{BE}, where $q=0$ and $h(t)\equiv e^{\la
t}$.

 \begin{theorem}\label{tq2ex}
 Let $f\in \elle\infty$, and let $u$ be a solution of \rife{q}, with $0\leq q<2$. Assume
that
\be\label{exp}
\lim\limits_{s\to +\infty}\frac{h(s)}{\int_0^s h(t)dt}=\la>0\,,\quad \hbox{for every $t\in
\r$, \ \ \ $\exists \lim\limits_{s\to +\infty} \frac{h(s+t)}{h(s)}\,:\,= e^{\la t}$.}
\ee
 Then we have:
 \be\label{ldir3}
\lim\limits_{\de \to 0}\de\, \frac{
\partial u}{
\partial
\nu(x)}(x-\de\nu(x))=\frac2{\la}\,,\qquad \lim\limits_{\de \to 0}\de\, \frac{\partial u}{
\partial \tau(x)}(x-\de\nu(x))=0
\ee and therefore 
\be\label{lglo3}
 \lim\limits_{\dom \to 0} \dom \nabla u(x)=\frac 2{\la} \nu.
\ee
 \end{theorem}
\Proof We use the same framework of the proof of Theorem \ref{tq2},  setting 
$$
v_\de=u(\de \xi)-\tilde F(\de), $$
where the function $\tilde F$ is defined by
\be\label{estesp'}\tilde F^{-1}(s)=
 \int^{+\infty}_s \frac{1}{[2\int_0^t h(\xi)d\xi]^{\frac12}}dt.
\ee
Indeed, as a consequence
of Keller-Osserman estimate and due to \rife{exp}, there holds
\be\label{estesp}
 |u(x)-\tilde F(\dom)|\leq C_0\quad \hbox{for any $x\in \Omega$: $\dom \leq d_0$, .}
 \ee 
 Observe that, since $\lim\limits_{s\to +\infty}\frac{h(s)}{\int_0^s h(t)dt}=\la>0$,
one can prove (as in \rife{dlog}) that $\tilde F'(t)t$ is bounded on $\r^+$ and 
\be\label{dert} \tilde F'(\de)\de\to - \frac2{\la}\qquad\hbox{as $\de\to0$.}
\ee Moreover 
the function $\tilde F$ is convex, so that we still have \rife{prof}, and then  again
\be\label{alw2}
 |u(\de\xi)- \tilde F(\de\,(\xi_1+\de^{1-2\sigma}))| \leq C_1\quad \hbox{for any $\xi\in
\frac1\de D_\de$.}
 \ee Reasoning as in the proof of Theorem \ref{tq2} we deduce that there exist   positive
constants $\ga_0$, $\ga_1$, $\sigma_0$ such that $$
\tilde F(\de(\xi_1+\de^{1-2\sigma}))-\tilde F(\de)\geq -\ga_0\log (\xi_1+\de^{1-2\sigma}) \qquad
\hbox{if $\xi_1\leq 1-\de^{1-2\sigma}$,} $$
and $$
\tilde F(\de(\xi_1+\de^{1-2\sigma}))-\tilde F(\de)\geq -\ga_1\log (\xi_1+\de^{1-2\sigma}) \qquad
\hbox{if $1<\xi_1< \frac{\sigma_0}\de-\de^{1-2\sigma}$,} $$
which together with \rife{alw2} imply 
\be\label{zer2} v_\de(\xi)\geq -\ga_0\log(\xi_1+\de^{1-2\sigma})- c_0\qquad \hbox{if
$\xi_1\leq 1-\de^{1-2\sigma}$,}
\ee and
\be\label{inf2} v_\de(\xi)\geq -\ga_1\log (\xi_1+\de^{1-2\sigma}) -c_1\qquad \hbox{if
$1<\xi_1< \frac{\sigma_0}\de-\de^{1-2\sigma} $.}
\ee Now the function $v_\de$ satisfies the equation 
$$
-\Delta v_\de +h(u(\de\xi)) \de^2+|\nabla v_\de|^q\,\de^{2-q}=\de^2 f(\de\xi)\qquad \xi
\in \frac1\de
 D_\de 
$$
and $v_\de$ is locally uniformly bounded. Since 
$$
\de^2 h(u(\de\xi))= \frac{h(v_\de+\tilde F(\de))}{h(\tilde F(\de))}\,
\frac{h(\tilde F(\de))}{\int_0^{\tilde F(\de)}h(s)ds}\, e^{2\log (\de
[\int_0^{\tilde F(\de)}h(s)ds]^{\frac12})}\,,
$$
as in the proof of Theorem \ref{tq2} we obtain, using \rife{exp} and \rife{dert}, 
that $ \de^2 h(u(\de\xi))$ is locally uniformly bounded and moreover
$$
\lim\limits_{\de \to 0}\de^2 h(u(\de\xi))= e^{\la v}\frac2\la\, 
$$
locally uniformly, where $v$ is the limit of  a subsequence (not relabeled) of $v_\de$. 
When $q>1$, local estimates of Bernstein's type (see e.g. \cite{LL}, \cite{Li} and the remark therein of the regularity of $f$), 
imply that any solution of \rife{q} satisfies, for a  constant $C>0$,
$$
|\nabla u(x)| \leq C \dom^{-\frac 1{q-1}}\,.
$$ 
In particular $v_\de$ verifies an equation of type
\be\label{trans}
-\Delta v_\de+ F_\de\cdot \nabla v_\de= g_\de\,,
\ee
where $g_\de$, $F_\de$ are  a function, and a field respectively, which are locally uniformly bounded. By elliptic estimates
we deduce that $\nabla v_\de$ is also locally uniformly bounded, and $v_\de$ is relatively compact in the $C^1_{loc}$-topology.  We have therefore
$$
\lim\limits_{\de \to 0}|\nabla v_\de|^q\, \de^{2-q}=0\,.
$$
When $0\leq q\leq  1$, $u\in L^{\infty}_{loc}(\Omega)\cap H^1_{loc}(\Omega)$ implies $\abs{\nabla u}^q\in L^{2/q}_{loc}(\Omega)$. Thus, by elliptic equations regularity theory and a  standard bootstraping argument, it follows that $\nabla u$ remains locally bounded   and the above limit holds true directly.
Thus, by replacing $g_{\de}$ by its expression and using also \rife{zer2}--\rife{inf2}, it turns out that $v$ is a solution of 
\be\label{espv}
\cases{ -\Delta v+\frac2\la\,e^{\la v}=0 & in $\rnp$, \cr \lim\limits_{\xi_1\to 0^+}v(\xi)
=+\infty\,,\cr}
\ee satisfying in addition that there exists $\ga$, $C>0$ such that for any $S>1$ we have
\be\label{addv} v(\xi)\geq -\ga\log S -C\qquad \hbox{for any $\xi$:  $\xi_1\leq S$.}
\ee 
When By Lemma \ref{sym} we conclude that $v= -\frac2\la\log\xi_1$, and this uniqueness result implies also that the whole sequence $v_\de$ is converging in $C^1_{loc}(\rnp)$. The convergence of
$\nabla v_\de$ to $\nabla v$ then yields \rife{ldir3} and \rife{lglo3}.
\qed
\vskip1em

\begin{remark} \rm  As a byproduct  of the scaling argument, from the convergence of
$v_\de=u(\de\xi)-\tilde F(\de)$ to $-\frac2\la\log\xi_1$, we obtained, setting $\xi=1$, that 
$$
 u(x)-\tilde F(\dom) \to 0\qquad\hbox{as $\dom\to 0$,} $$
where $\tilde F$ is defined in \rife{estesp'}. In case $q=0$ we recover a result of  \cite{LMK}.
\end{remark}

Finally, we consider the  case that $h$ has  a power--type asymptotic rescaling at
infinity: 
we extend then some results  proved in \cite{BaM2} for the case $q=0$.
\vskip1em

 \begin{theorem}\label{tq2pow}
 Let $f\in \elle\infty$ and $u$ be a solution of \rife{q}, with $0\leq q<2$.  
 \begin{itemize}
 \item[(i)]
 Assume that
\be\label{surc}
\lim\limits_{s\to +\infty}\frac{h(s)^{\frac2q}}{\int_0^s h(t)dt}=+\infty\,,
\ee 
\end {itemize}
and
\be\label{scapow}
\quad \hbox{for every $t\in
\r^+$ \ \ \ $\exists \lim\limits_{s\to +\infty} \frac{h(st)}{h(s)}\,:\,=  t^\alpha$, with
$\alpha>1$.}
\ee
 Then we have:
 \be\label{ldir4}
\lim\limits_{\de \to 0}\frac1{\tilde F'(\de)}\, \frac{
\partial u}{
\partial
\nu(x)}(x-\de\nu(x))=  1
\,,\qquad \lim\limits_{\de \to 0}\frac1{\tilde F'(\de)}\, \frac{
\partial u}{
\partial \tau(x)}(x-\de\nu(x))=0
\ee
where $\tilde F^{-1}(s)$ is defined in \rife {estesp'}, and in
particular 
\be\label{lglo4}
 \lim\limits_{\dom \to 0}  \frac{\nabla u(x)}{\tilde F'(\dom)}= \nu.
 \ee
  \begin{itemize}
 \item[(ii)] Assume that $q>1$ and 
 \be\label{sousc}
\lim\limits_{s\to +\infty}\frac{h(s)^{\frac2q}}{\int_0^s h(t)dt}=l\,,
\ee 
\end {itemize}
for some $l\geq 0$, and let $a>0$ be  such that $\frac a{2-q}-a^{\frac q2}= (\frac  {2-q} 2l)^{\frac
q{2-q}}$. Then 
 \be\label{ldir5}
\lim\limits_{\de \to 0} \de^{\frac1{q-1}}\, \frac{
\partial u}{
\partial
\nu(x)}(x-\de\nu(x))=  b_q
\,,\qquad \lim\limits_{\de \to 0}\de^{\frac1{q-1}}\, \frac{
\partial u}{
\partial \tau(x)}(x-\de\nu(x))=0
\ee 
where
$b_q=\left(\frac1a\right)^{\frac{2-q}{2(q-1)}}\left(\frac{2-q}{q-1}\right)^{\frac1{q-1}}
$, and then 
\be\label{lglo5}
\lim\limits_{\dom \to 0}  \dom^{\frac1{q-1}}\nabla u(x)= b_q\nu.
\ee
 \end{theorem}

\begin{remark}\rm As pointed out in Remark \ref{tre}, 
the existence of the limit in \rife{scapow} automatically implies 
that this limit is a power function.
\end{remark}
\Proof
(i) Under assumption \rife{surc},  we can apply the results in \cite{BG} and use that
\be\label{rapp}
\lim\limits_{\dom \to 0} \frac{u(x)}{\tilde F(\dom)}=1,
\ee 
In other words, the behaviour of $u$ is determined by the Keller--Osserman estimate in
this case. Let us now use the framework of Theorem \ref{tq2}, introducing the system of
coordinates $(\eta_1, \ldots,\eta_N)$ whose $\eta_1$--axis is the inner normal direction.
Define $O_\de=(\de,\ldots,0)$ and the domain $$
\tilde D_\de= B(O_\de,\de^{1-\sigma})\cap B(P_0, d_0-\de)\,,\qquad \sigma\in (0,\frac12).
$$
Again we have that $\tilde D_\de$ converges to the half space $\{\xi\,:\,\xi_1>0\}$. Now
we set $\xi=\frac{\eta-O_\de}\de$ and we introduce the blown-up function 
$$
v_\de=\frac {u(\de\xi+O_\de)}{\tilde F(\de)}. $$
This time let us choose $d_0$ such that $\dom<d_0$ implies $|\frac{u(x)}{\tilde F(\dom )}-1|\leq
\vep_0$; thanks to \rife{Tay} it follows 
$$
(1-\vep_0)\tilde F(\de(\xi_1+1)+ O(\de^{2-2\sigma}))\leq u(\de\xi+O_\de)\leq
(1+\vep_0)\tilde F(\de(\xi_1+1)+ O(\de^{2-2\sigma}))\,. $$ 
In particular we deduce that  $0\leq v_\de\leq (1+\vep_0)$, i.e. $v_\de$ is uniformly
bounded and satisfies $$
-\Delta v_\de +\frac{h(u(\de\xi+O_\de))\de^2}{\tilde F(\de)}+\tilde F(\de)^{q-1}\,|\nabla v_\de|^q
\de^{2-q}= f(\de\xi+O_\de)\frac{\de^2}{\tilde F(\de)}\,. $$
Note that \rife{scapow} implies
\be\label{meg}
\lim\limits_{t\to +\infty} \tilde F^{-1}(t)\sqrt{\frac {h(t)}t}=\lim\limits_{t\to +\infty}
\int_1^{+\infty} \frac1{\sqrt{2\int_0^s\frac{h(t\xi)}{h(t)}d\xi
}}ds=\frac{\sqrt{2(\alpha+1)}}{\alpha-1}
\ee so that $$
\lim\limits_{\de\to 0} \frac{h(\tilde F(\de))\de^2}{\tilde F(\de)}= \frac{ 2(\alpha+1)}{(\alpha-1)^2}.
$$
Set $c_\alpha= \frac{ 2(\alpha+1)}{(\alpha-1)^2}$; then we have, using that $v_\de$ (up to
subsequences) converges,  locally uniformly, to a function $v$, and $\frac{h(st)}{h(s)}$
converges to $t^\alpha$ locally uniformly in $\r$, 
$$
\frac{h(u(\de\xi+O_\de))\de^2}{\tilde F(\de)}=\frac{h(v_\de
\tilde F(\de))}{h(\tilde F(\de))}\,\,\frac{h(\tilde F(\de))\de^2}{\tilde F(\de)}\to 
c_\alpha\, v^\alpha\,. 
$$
As in the previous theorem, we can use the local estimates on $\nabla u$ for solutions of  \rife{q}, in order to get
$$
|\nabla u(x)|\leq C\dom^{-\frac1{q-1}}
$$ 
when $q>1$ for some constant $C>0$. This implies that $\tilde F(\de)^{q-1}\,|\nabla v_\de|^{q-1}
\de^{2-q}$ is locally uniformly bounded. Hence $v_\de$ satisfies an equation like \rife{trans} with $g_\de$ and $F_\de$ locally bounded. We deduce with a  simple bootstrap argument and elliptic regularity  that $v_\de$
is relatively compact in the $C^1_{loc}$-topology. Moreover assumption \rife{surc} implies that 
$$
\lim\limits_{t\to +\infty} 
\frac 1{t^{\frac2{2-q}}}\,\int_0^t h(s)ds=+\infty\,, $$
which in turn gives that 
$$
\lim\limits_{\de\to 0}  \tilde F(\de)^{q-1}\, \de^{2-q}=0\,. 
$$
Therefore we conclude that $$
\tilde F(\de)^{q-1} |\nabla v_\de|^q \de^{2-q}\quad\mathop{\to}^{\de\to 0} \quad 0\,. $$
When $0\leq q\leq 1$, $\nabla u$ remains locally bounded and the same conclusion holds.
In both case we conclude that the function $v$ satisfies, in the limit, the equation
\be\label{vpow} -\Delta v + c_\alpha v^\alpha=0 \qquad \hbox{in $\rnp$,}
\ee and it is uniformly bounded. 

By \rife{scapow} and the dominated converge theorem, 
\be\label{ehh}
\lim\limits_{\xi\to+\infty}\frac{\int_0^\xi h(s)ds}{\xi\,h(\xi)}=
\lim\limits_{\xi\to+\infty}\frac{\int_0^1 h(\xi s)ds}{h(\xi)}
=\int_0^1\lim\limits_{\xi\to+\infty}\frac{h(\xi s)}{h(\xi)} ds
= \frac 1{\al+1},
\ee 
then, using  \rife{meg}, there holds
\be\label{hom}
\eqalign{ &
\lim\limits_{t\to 0} \frac{t\tilde F'(t)}{\tilde F(t)}= 
\lim\limits_{\xi\to +\infty}-
\frac{\tilde F^{-1}(\xi)\sqrt{2\int_0^\xi h(s)ds}}\xi=\cr & =\lim\limits_{\xi\to
+\infty}-\tilde F^{-1}(\xi)\, \sqrt{\frac{h(\xi)}\xi} \sqrt{\frac{2\int_0^\xi
h(s)ds}{\xi\,h(\xi)}}= -\frac{\sqrt{2(\alpha+1)}}{\alpha-1}\,\sqrt{\frac2{\alpha+1}}  =-
\frac2{\al-1}\,.\cr}
\ee Moreover, the function $\frac{\tilde F'(\xi)}{\tilde F(\xi)}$ is increasing, so that for any
$\la>1$, $$
1\geq \frac{\tilde F(\la \de)}{\tilde F(\de)}=\exp\left[\log(\tilde F(\la\de))-\log(\tilde F(\de))\right]\geq
\exp\left[\frac{\tilde F'(\de)\de}{\tilde F(\de)}(\la-1)\right]\geq \vep_0>0\,. $$
Thus, for any $\la>1$ the sequence $\frac{\tilde F(\la\,\de)}{\tilde F(\de)}$ is bounded, strictly
positive, and satisfies, in  view of \rife{hom} and \rife{ehh}, $$
\frac{\tilde F(\la\,\de)}{\tilde F(\de)}=\frac{\tilde F'(\la\,\de)}{\tilde F'(\de)}\,\la(1+o(1))=
 \sqrt{\frac{\int_0^{\tilde F(\la\de)}h(s)ds}{\int_0^{\tilde F(\de)}h(s)ds}}\,\la(1+o(1))=
 \sqrt{\frac{ \tilde F(\la\de) }{\tilde F(\de)}\,\frac{h(\tilde F(\la\,\de))}{h(\tilde F(\de))}}\,\la(1+o(1))\,.
 $$
Using  \rife{scapow} we deduce 
\be\label{hom2}
\exists \,\, \lim\limits_{\de \to 0}\frac{\tilde F(\la\,\de)}{\tilde F(\de)}=\la^{-\frac2{\al-1}}\,.
\ee Then   we have 
\be\label{basso}
\frac{\tilde F(\de(\xi_1+1+ O(\de^{1-2\sigma}))}{\tilde F(\de)}\geq 
 (1+\xi_1)^{-\frac2{\al-1}}+o(1) \quad \hbox{as $\de\to 0$.}
\ee Since we have $$
\frac{u(\de\xi+O_\de)}{\tilde F(\de(\xi_1+1)+ O(\de^{2-2\sigma}))}\,\frac{\tilde F(\de(\xi_1+1)+
O(\de^{2-2\sigma}))}{\tilde F(\de)}\leq v_\de\leq 
\frac{u(\de\xi+O_\de)}{\tilde F(\de(\xi_1+1)+ O(\de^{2-2\sigma}))}\, $$
from \rife{rapp} (recall that $\eta=\de\xi+O_\de$ and $\dist(\eta,\partial
\Omega)$ is estimated in  \rife{Tay}) and
\rife{basso} we obtain: 
$$
 (1+\xi_1)^{-\frac2{\al-1}}+o(1)\leq v_\de(\xi)\leq 1+o(1)
\quad\hbox{as $\de\to0$,} 
$$
and we conclude that 
$$
\lim\limits_{\xi_1\to 0^+} v(\xi_1)=1. $$
Together with \rife{vpow} this implies that 
$v=(1+\xi_1)^{-\frac2{\alpha-1}}$. The $C^1_{loc}$ convergence of $v_\de$ gives then  
$$
\nabla v_\de(\xi)\to -\frac2{\alpha-1}(1+\xi_1)^{-\frac{\alpha+1}{\al-1}}(1,
0,\ldots,0)\,. $$
Now recall that $\nabla_\xi u(\de\xi+0_\de)=\
\frac{\tilde F(\de)}\de \nabla v_\de(\xi)$, hence using \rife{hom}--\rife{hom2} we get $$
\frac{\nabla_\xi u(\de\xi+0_\de)}{\tilde F'(\de(1+\xi_1))}\to (1, 0,\ldots,0)\,, $$
which gives \rife{ldir4} and \rife{lglo4}.

(ii) Using  \rife{sousc}, we have from \cite{BG} and \cite{G}:
\be\label{rapp2}
\lim\limits_{\dom \to 0} \frac{u(x)}{c_q\dom^{-\frac{2-q}{q-1}}}=1\,,
\ee where $c_q = \left( \frac{2-q}{(q-1)\sqrt a}\right)^{\frac{2-q}{q-1}}$. With the same
notations as above we set $$
v_\de(\xi)= \frac{u(\de\xi+O_\de)}{c_q \de^{-\frac{2-q}{q-1}}}. $$
As before, we deduce that $v_\de$  is uniformly bounded, and satisfies $$
-\Delta v_\de +\frac{h(u(\de\xi+O_\de))\de^{\frac q{q-1}}}{c_q  }+c_q^{q-1}\,|\nabla
v_\de|^q = f(\de\xi+O_\de)\frac{\de^{\frac q{q-1}}}{c_q }\,. $$
 Now assumption \rife{sousc} implies
\be\label{ehhh}
\lim\limits_{t\to +\infty}
\frac1{t^{\frac2{2-q}}}\,\int_0^t h(s)ds=\left(\frac{2-q}2\,l^{\frac q2}\right)^{\frac
2{2-q}}.
\ee 
Noticing that
$$
h(u(\de\xi+O_\de))\de^{\frac q{q-1}}=
\frac{h(u(\de\xi+O_\de))}{\left(\int_0^{u(\de\xi+O_\de)}h(s)ds\right)^{\frac q2}}\,
\frac{\left(\int_0^{u(\de\xi+O_\de)}h(s)ds\right)^{\frac q2}}{u(\de\xi+O_\de)^{\frac
q{2-q}}}
 (c_q\, v_\de)^{\frac q{2-q}} \,, $$
and using \rife{ehhh} and assumption \rife{sousc},  we get 
$$
h(u(\de\xi+O_\de))\de^{\frac q{q-1}}\mathop{\to}^{\de\to 0}(\frac{2-q}2\,l)^{\frac
q{2-q}}(c_q\, v )^{\frac q{2-q}}\,\,. $$
Therefore passing to the limit as $\gd\to 0$, we conclude that $v$ solves 
\be\label{lav} -\Delta v + (\frac{2-q}2\,l)^{\frac q{2-q}}\,c_q^{\frac q{2-q}-1}\,
v^{\frac q{2-q}}+c_q^{q-1}\,|\nabla v|^q =0.
\ee Similarly as for (i), thanks to \rife{rapp2} we also obtain that
\be\label{liv}
\lim\limits_{\xi_1\to 0^+} v(\xi)=1\,,\quad \lim\limits_{\xi_1\to+\infty}v(\xi)=0\,.
\ee Recalling the value of $c_q$ and the definition of $a$ in \rife{sousc}, 
one can check that the function $(\xi_1+1)^{-\frac{2-q}{q-1}}$ is  a solution of 
\rife{lav}--\rife{liv}. On the other hand, for any $\alpha\geq 0$, $\beta>0$, the problem 
\be\label{lpr}
\cases{-\Delta z+\al z^{\frac q{2-q}}+\beta\,|\nabla z|^q =0&in $\rnp$,\cr
z_{\mathop{|}_{\xi_1=0}}=1\,,\quad \lim\limits_{\xi_1\to +\infty} z(\xi)=0\,&\cr}
\ee
admits one and only one positive solution:  see Theorem \ref{sym2} below for  a   more general
result of this type.

Having an explicit solution of \rife{lav}--\rife{liv}, we conclude  that
$v=(\xi_1+1)^{-\frac{2-q}{q-1}}$. The uniqueness of this limit yields the convergence of
the whole sequence $v_\de$, in particular we get that $\nabla v_\de(\xi)$ converges to 
$\nabla v(\xi)$ locally uniformly. Setting $\xi=(1, 0, \ldots, 0)$, we obtain relations
\rife{ldir5}--\rife{lglo5}.
\qed

\begin{remark} \rm  The  result of Theorem \ref{tq2ex} still holds if one  relax the
assumptions on the right hand side: for the case (i), it is enough to  require that
$\lim\limits_{\dom \to 0}
\frac{f(x)}{h(\tilde F(\dom))}=0$, where $\tilde F$ is defined through \rife{estesp'}. Note that if
$h(s)=|s|^{\beta-1}s$ ($\beta>1$), this means $\lim\limits_{\dom \to 0}\dom^{\frac
{2\beta}{\beta-1}} f(x)=0$.

In case (ii), it would  be enough to have $\lim\limits_{\dom \to 0}\dom^{\frac q{q-1}}
f(x)=0$; in fact, this corresponds to the case $h(s)=s^\beta$, with $\beta\leq \frac
q{2-q}$.
\end{remark}

\begin{remark} \rm  
In case $h(u)=\la\, u$, the (unique) solution of \rife{equa} is the value function of an associated suitable stochastic control problem with state constraint, which is described in \cite{LL}. In that context, the field $-q|\nabla u|^{q-2}\nabla u$ is exactly the optimal feedback control, whose role is to keep the process to stay inside $\Omega$ (minimizing a  certain cost functional). Our results (Theorem \ref{tq2} and Theorem \ref{tq2pow}) prove the precise asymptotics for the control, i.e. 
$-q|\nabla u(x)|^{q-2}\nabla u(x) \sim -\frac {q'}{\dom}\,\nu(x)$ as $\dom \to 0$. 

\end{remark}

\section{On the  uniqueness of explosive solutions in case of concavity}

In this section we give a   uniqueness  result for solutions of
\be\label{pb}
 \cases{-\Delta u+h(u) +|\nabla u|^q = f& in $\Omega$, \cr \lim\limits_{\dom\to 0}u(x)
=+\infty\,,&\cr}
\ee which applies to the case that  $h(s)$ is  concave. 
We restrict ourselves to $q>1$, which is the significant case. Our basic criterion for
uniqueness is the following.

\begin{theorem}\label{criter} Let $\Omega$ be  a bounded domain and $f\in \elle\infty$.
Assume $1<q\leq 2$, and that $h$ is a  continuous increasing function satisfying the
following assumption:
\be\label{conh}
\begin{array}{rl}&
\hbox{$\exists$ a  positive, continuous function $m(s)$, and constants $c_0$, $\vep_0>0$
such that} \\[2mm]

& h((1+\vep)a+\vep b)-(1+\vep)h(a)\geq \vep \,b\,m(a)-c_0\,\vep(1+a)\,,\\[2mm]&
\quad \hbox{$\forall a\in \r$, 
$\vep\in (0,\vep_0)$, $0\leq b\leq \frac1{\vep_0}$.}
\end{array}
\ee If $u_1$, $u_2$ are two solutions of \rife{pb} such that
\be\label{apri}
\lim\limits_{\dom\to 0}\frac{u_1}{u_2}=1\,,\qquad \lim\limits_{\dom\to 0}\frac{|\nabla
u_i|^q}{u_i}=+\infty\quad\forall i=1,2,
\ee then $u_1=u_2$.
\end{theorem}
\Proof
We set $A(v)=-\Delta v+h(v)+|\nabla v|^q$. Define   $u_2^\vep=(1+\vep)u_2+\vep T$, where
$T$ is  a positive constant  to be chosen later. Then 
$$
A(\ude)= \,h((1+\vep)u_2+\vep T)-(1+\vep)h(u_2)+(1+\vep) ((1+\vep)^{q-1}-1)|\nabla u_2|^q
+(1+\vep)f\,, $$
and using \rife{conh} and that $f\in \elle\infty$
\be\label{pred}
\eqalign{ & A(\ude)\geq   
f+\vep[ m(u_2)T+(q-1)|\nabla u_2|^q-c_1(1+u_2)]\,.\cr}
\ee By assumption \rife{apri}, there exists a positive, bounded, compactly supported
function $\psi(x)$ such that $$
(q-1) |\nabla u_2|^q-c_1(1+u_2) \geq -\psi(x)\,. $$
If $K\subset \Omega$ is a compact set containing the support of $\psi$, we have that $u_2$ is
bounded on $K$ and since $m(s)$ is positive we have  $\inf\limits_K m(u_2) >0$. Setting
$T=\frac {\|\psi\|_{\infty}}{\inf\limits_K m(u_2)}$ then implies $$
A(\ude)\geq f=A(u_1)\qquad \hbox{in $\Omega$.} $$
Moreover since $\frac{u_1}{u_2}\to 1$ as $\dom\to 0$, we have that $u_1-\ude\leq 0$ near $
\partial \Omega$. Inside $\Omega$, 
we use that $h$ is increasing  to deduce that $u_1-\ude\leq 0$ on any maximum point, so
that we can conclude that $$
 u_1\leq (1+\vep)u_2+\vep T     \qquad \hbox{in $\Omega$.} $$
Letting $\vep \to 0$ we get $u_1\leq u_2$. Interchanging the roles of $u_1$, $u_2$, we
conclude that $u_1=u_2$
\qed

Let us make some comments and remarks about the previous result:
\begin{itemize}

\item[1)] Assumption \rife{conh} is satisfied if $h(s)=h_1(s)+h_2(s)$, where $h_1$ is  a
nondecreasing convex function and $h_2$ is an increasing concave function. Indeed, one
has, taking into account the sublinear behaviour of the concave part, 
$$
\begin {array}{l} 
h((1+\vep)a+\vep b)-(1+\vep)h(a) \geq -\vep h_1(0)+ h_2((1+\vep)a+\vep
b)-h_2((1+\vep)a)\\ 
\phantom{A((1+\vep)u_2+\vep T)...................} \qquad\qquad\qquad
+h_2((1+\vep)a)-h_2(a)- \vep h_2(a)\\
\phantom{A((1+\vep)u_2+\vep T)...................}\geq \vep
\,b\,m(a)-c_0\,\vep(1+a)\,,
\end {array} $$
with, for instance,   $m(a)=h_2'(2a+1)$ if $h$ is differentiable, or
$m(a)=h_2(2a+1)-h_2(2a)$ otherwise. 

\item[2)] As remarked above, the previous result 
is meant to  apply to the case that $h$ is the sum of  a convex function and an {\sl 
increasing} concave function. On the other hand, we recall that in case $h$ is purely convex the
uniqueness of solutions has been proved in previous papers (see e.g. \cite{LL}),
essentially using the following  standard argument: if $\frac{u_1}{u_2}\to 1$ as $\dom\to
0$, then it is enough to take $T>-h^{-1}(m)$ where $m=\inf_{\Omega} f$, in order to have
$$\begin {array}{l}
A((1+\vep)u_2+\vep T)\geq \,h((1+\vep)u_2+\vep T)-(1+\vep)h(u_2)\\
\phantom{A((1+\vep)u_2+\vep T)...................}
+(1+\vep) ((1+\vep)^{q-1}-1)|\nabla u_2|^q +(1+\vep)f\\
\phantom{A((1+\vep)u_2+\vep T)}
\geq f+\vep[  f-h(-T)]\\\phantom{A((1+\vep)u_2+\vep T)}
\geq f= A(u_1), 
\end {array}$$
which yields $u_1\leq u_2$ for any $u_1$, $u_2$ large solutions such that
$\frac{u_1}{u_2}\to 1$ as $\dom\to 0$. Note that in this case one does not need to have
any information with respect to the gradients.

 \item[3)] Assumption \rife{apri} is not really restrictive, and is  certainly satisfied
in smooth domains $\Omega$ and in almost all significant situations. Indeed, this is a
consequence of the results on the asymptotic behaviour of $u$ and $\nabla u$ which are
given in Section 2, so that in particular \rife{apri}  is   verified for all  the
situations considered in  Theorem  \ref{tq2}, Theorem \ref{tq2ex} and  Theorem
\ref{tq2pow}, which deal with possibly power or exponential growths of $h$ at infinity.  
 
In particular, this applies to the case that  $h$ is concave (which implies assumption
(h1) in Theorem \ref{tq2} and assumption \rife{sousc} in Theorem \ref{tq2pow}), hence  
condition \rife{apri} follows  from Section 2 and \rife{conh} also holds true. We get then
the  following  corollary.

\begin{corollary} Let $\Omega$ be a smooth domain and $f\in \elle\infty$. If $h$ is
increasing and concave, for any $q>1$ problem \rife{pb} has a unique solution.
\end{corollary}
 
On the other hand, note that   for possibly larger growths of $h$ than considered in
Section 2,  more precisely when $$
 \hbox{either $q=2$  and $ \lim\limits_{s\to +\infty} \frac{h(s)e^{-2s}}{\int_0^s
h(t)e^{-2t}dt}=+\infty$, or $q<2$ and $ \lim\limits_{s\to +\infty} \frac{h(s) }{\int_0^s
h(t) dt}=+\infty$,} $$
 uniqueness of large solutions follows easily since one can prove directly that
$u_1(x)-u_2(x)\to 0$ as $\dom \to 0$ for any two solutions $u_1$, $u_2$. Therefore the
problem of uniqueness is really significant when $h$ satisfies growth conditions of the same
kind as in Section 2.
 
 \end{itemize}

\section{Appendix: On some  symmetry results in the half space}

In the proof of Theorem \ref{tq2pow} we have used  a uniqueness result for solutions of \rife{lpr}.
Here we give  a self--contained proof of a even more general  result on the uniqueness, or symmetry, of nonnegative solutions of such type of problems in the half space, without conditions at infinity. More precisely, consider the problem
\be\label{gen}
\cases{-\Delta z+\al z^{p}+\beta\,|\nabla z|^q =0&in $\rnp\,:\,=\{\xi=(\xi_1,\xi')\in \rn\,:\,\xi_1>0\}$,\cr
\noalign{\medskip}
z\geq 0 & in $\rnp$,\cr
\noalign{\medskip}
\lim\limits_{\xi_1\to 0^+} z(\xi_1,\xi')=M&locally uniformly with respect to $\xi'$\,,\cr}
\ee
where $0\leq M\leq \infty$,  $\beta$, $p>0$, and $\alpha\geq 0$.

Next we prove that the solutions of \rife{gen} are one--dimensional, and in particular unique if $\alpha >0$.
\vskip1em

\begin{theorem}\label{sym2}
Let $1<q\leq 2$,  $\al\geq 0$, $p>0$ and $\beta>0$. Let also $M\in [0,+\infty]$. Then  
\vskip0,5em
(i) if $\al>0$ problem \rife{gen} admits  a unique solution $z$, and $z= z(\xi_1)$ 
\vskip0,5em
(ii) if  $\al=0$ any solution   of \rife{gen} is a function of the only variable $\xi_1$. In particular, 

(a) if $q=2$ then  (necessarily)   $M<\infty$ and  $z\equiv M$. 

(b) if $q<2$ and $M<\infty$ then either $z\equiv M$  or there exists $l\in [0,M)$ such that  
$$
z(\xi)= l+\int_{\xi_1}^{+\infty} [(q-1)s+c_{M,l}]^{-\frac1{q-1}}ds\,,$$
where $c_{M,l}$ is uniquely determined by the implicit relation
$$\int_{0}^{+\infty} [(q-1)s+c_{M,l}]^{-\frac1{q-1}}ds=M-l\,,
$$
while if $M=+\infty$ then there exists $l\in [0, +\infty)$ such that $$z(\xi)= l+ \int_{\xi_1}^{+\infty} [(q-1)s ]^{-\frac1{q-1}}ds.$$ 
\end{theorem}
\Proof 
  (i) Let $\al>0$.  First of all, as in Lemma \ref{sym},  consider the radial solutions $\omega_R$ of 
$$
\cases{-\Delta \omega_R+\al \omega_R^{p}+\beta\,|\nabla \omega_R|^q =0&in $B_R(0)$,\cr
\lim\limits_{\rho\uparrow R} \omega_R(\rho)=+\infty\,, &\cr}
$$ 
and the sequence $\{\omega_R(\xi-\xi_R)\}_{R>0}$, where $\xi_R=(R,0)$. 
Note that this sequence exists since $\al>0$ and $q>1$. By local estimates we have that $\omega_R(\cdot-\xi_R)$ is locally bounded and moreover it is a decreasing sequence converging towards  a function $\omega_{\infty}(\xi_1)$ which is the unique solution of 
\be\label{1di}
\cases{- \omega_\infty''+\al \omega_\infty^{p}+\beta\,|  \omega_\infty'|^q =0 \quad \hbox{in $(0,+\infty)$}&\cr
\omega_\infty(0)=+\infty\,.&\cr}
\ee
Indeed, $\omega_\infty$ is a positive, decreasing convex function and converges to zero as $\xi_1$ tends to infinity. 
Since any solution $z$ of \rife{gen} is below $\omega_R(\xi-\xi_R)$ on $B_R(R,0)$, we deduce in the limit that
\be\label{hain}
z\leq \omega_{\infty}(\xi_1)\,.
\ee
In particular, $z$ tends to zero as $\xi_1$ tends to infinity. Now, for $R,\,S>0$, consider the radial solutions ${\underline \omega}_{R,S}(\rho)$
of 
$$
\cases{-\Delta{\underline \omega}_{R,S}+\al {\underline \omega}_{R,S}^{p}+\beta\,|\nabla {\underline \omega}_{R,S}|^q =0&in $B_{R+S}(0)\setminus B_R(0)$,\cr
\lim\limits_{\rho\downarrow R} {\underline \omega}_{R,S}=M\,, \qquad {\underline \omega}_{R,S}(R+S)=0&\cr}
$$
and the sequence $\{ {\underline \omega}_{R,S}(\xi - \eta_R)\}_{R,S}$, where $\eta_R=(-R,0)$. 
It can be easily checked that, since ${\underline \omega}_{R,S}$ is positive and decreasing
with respect to $\rho$,  the sequence
$\{ {\underline \omega}_{R,S}(\cdot - \eta_R)\}_{R,S}$ is increasing respect to $R$ and $S$. Letting successively $R\to\infty$ and $S\to\infty$, its limit $\omega_M$ is 
 a one--dimensional solution of \rife{gen}. By comparison we
have that
$\{ {\underline \omega}_{R,S}(\xi - \eta_R)\}_R\leq z(x) $, for any solution $z$
 of \rife{gen}, hence we get in the limit
\be\label{mu}
\omega_M(\xi_1)\leq z(x)\,\quad\forall \xi_{1}>0.
\ee
Now, since $\al>0$ the one dimensional solution of  \rife{gen} is unique; thus if $M=+\infty$, we have obtained that $z\equiv \omega_\infty(\xi_1)$. 

If $M<\infty$, we need a sharper upper bound for $z(x)$. To this purpose, let $t\in
(0,1)$;  we write $\xi=(\xi_1,\xi')$ and denote $B_R^{N-1}= \{|\xi'|<R\}\subset \bf
R^{N-1}$. We are going to construct a  supersolution in the cylinder $(0,L)\times
B^{N-1}_R$. 

\noindent Let   $\vfi_{t,L}(\xi_1)$ be the solution of the one-dimensional  problem
$$
\cases{
-\vfi_{t,L}''+ t^{p-1}\al \vfi_{t,L}^p+ t \beta |\vfi_{t,L}'|^q=0\quad \hbox{in
$(0,L)$},&\cr
\vfi_{t,L}(0)= \frac Mt\,\quad \vfi_{t,L}(L)=\frac1t\omega_\infty(L)\,,&\cr}
$$
where $\omega_\infty$ is the solution defined in \rife{1di}.
We also  set 
$$
f_t(s)=\cases{(1-t^2)^{\frac{p-1}2} s^p &if $p\geq 1$,\cr 
\noalign{\medskip}
\frac{
\left(\omega_\infty(L)+\sqrt{1-t^2} \,s\right)^p-\omega^p_\infty(L)}{\sqrt {1-t^2}}&if
$0<p<1$.\cr}
$$
Now consider the function
$\psi_{t,R}(\xi')$  solution of 
$$
\cases{-\Delta \psi_{t,R}+\al f_t(\psi_{t,R})+\beta\sqrt{1-t^2}
\,|\nabla \psi_{t,R}|^q =0&in $B^{N-1}_R\subset \bf R^{N-1}$,\cr
\lim\limits_{|\xi'|\uparrow  R}\psi_{t,R}=+\infty\,. &\cr}
$$
Note that such a function exists since $q>1$ and $f_t(s)$ is an increasing unbounded
function (in fact, $f_t(s)$ behaves like $(1-t^2)^{\frac{p-1}2} s^p$ for $s$ large). 

Define now $\bar z (\xi_{1},\xi')= t\vfi_{t,L}(\xi_{1})+\sqrt{1-t^2}\,\psi_{t,R}(\xi')$, we claim that $\bar z$ is a
supersolution. Indeed, using that $t\vfi_{t,L}\geq \omega_\infty(L)$ ($L$ is meant to be large enough so that $\omega_\infty(L)<M$), we have
$$
\bar z^p\geq t^p \vfi_{t,L}^p+ \sqrt{1-t^2}\, f_t(\psi_{t,R}).
$$
Moreover
$$
|\nabla \bar z|^q= \left( t^2 | \vfi_{t,L}' |^2+ (1-t^2)|\nabla
\psi_{t,R}|^2\right)^{\frac  q2}\geq  t^2 |  \vfi_{t,L}'|^q + (1-t^2) |\nabla
\psi_{t,R}|^q
$$
by concavity since $q\leq 2$, so that
$$
\eqalign{&
-\Delta \bar z+\al {\bar z}^{p}+\beta\,|\nabla \bar z|^q \geq 
t[ -\vfi_{t,R}''+ t^{p-1}\al \vfi_{t,L}^p+ t \beta |\vfi_{t,L}'|^q] \cr & \quad 
+ \sqrt{1-t^2} [ -\Delta \psi_{t,R}+\al
f_t(\psi_{t,R})+\beta\sqrt{1-t^2}\,|\nabla
\psi_{t,R}|^q]=0.\cr}
$$
Thus $\bar z$ is a supersolution of the equation in the cylinder  $(0,L)\times
B^{N-1}_R$. Moreover, since $\psi_{t,R}$ blows up at the boundary and is positive, and
using
\rife{hain}, we have that $z(x)\leq \bar z(x)$ on the boundary of the cylinder. By the
comparison principle we deduce that
$$
z(\xi)\leq t\vfi_{t,L}(\xi_1)+\sqrt{1-t^2}\,\psi_{t,R}(\xi')\qquad \hbox{in $(0,L)\times
B^{N-1}_R$.}
$$
Now let $R$ go to infinity, and use that $\psi_{t,R}$ converges to zero 
(as a consequence of the local estimates which depend on the distance to the boundary); we
obtain that
$$
z(\xi)\leq t\vfi_{t,L}(\xi_1)\,,
$$
and then, letting $L$ go to infinity, 
$$
z(\xi)\leq t\vfi_{t}(\xi_1)
$$
where $\vfi_t$ solves the problem 
$$
-\vfi_{t }''+ t^{p-1}\al \vfi_{t }^p+ t \beta |\vfi_{t }'|^q=0\quad 
\hbox{in $(0,+\infty)$},\quad \vfi_t(0)=\frac Mt.
$$
As $t$ tends to $1$, clearly $\vfi_t$ converges to the unique one--dimensional solution 
of \rife{gen}, which we called $\omega_M(\xi_1)$. Therefore $z\leq \omega_M(\xi_1)$, which
together with \rife{mu} gives the claimed result.
\vskip0,5em
(ii) Let now $\alpha=0$. Up to multiplying $z$ by a constant, we can assume that
$\beta=1$. We consider first the case $q<2$.

First observe that, since
$z$ is a solution in
$B_{\xi_1}(\xi_1,\xi')$, by the local estimates on $\nabla z$ (see e.g. \cite{LL}, \cite{Li}) we have
\be\label{gralo}
|\nabla z(\xi_1,\xi')|\leq C \xi_1^{-\frac1{q-1}}\qquad \forall (\xi_1,\xi')\in \rnp.
\ee
In particular,  we have
\be\label{cau}
|z(\eta_1,\xi')-z(\xi_1,\xi')|\leq  C\int_{\xi_1}^{\eta_1} t^{-{\frac1{q-1}}}dt
\ee
and since $\frac1{q-1}>1$ we deduce that $z(\xi_1,\xi')$ has a finite limit as $\xi_1$ goes to
infinity, and due to \rife{gralo} this limit does not depend on $\xi'$. Thus we set
$$
l\,:\,=\lim\limits_{\xi_1\to +\infty}u(\xi_1,\xi')\,.
$$
Using again \rife{cau} we also deduce the estimate:
\be\label{lest}
l- C\, \xi_1^{-\frac{2-q}{q-1}}\leq z(\xi_1,\xi') \leq l+ C\,\xi_1^{-\frac{2-q}{q-1}}
\qquad \forall (\xi_1,\xi')\in \rnp.
\ee
Our goal is now to prove that $z(\xi)=\omega_l(\xi_1)$, which is the unique solution of
$$
\omega_l''=|\omega_l '|^q\quad 
\hbox{in $(0,+\infty)$},\quad \omega_l(0)= M\,,\qquad \lim\limits_{\xi_1\to
+\infty}\omega_l(\xi_1)=l.
$$
In order to prove that $z\leq \omega_l$, let $t\in (0,1)$, $C\in \bf R$, and consider the problem on
$\bf R^{N-1}$: 
\be\label{erg} 
\cases{-\Delta \psi_{t,R}+\sqrt{1-t^2}\,|\nabla \psi_{t,R}|^q+ C=0& in $B_R^{N-1}\subset
\bf R^{N-1}$,\cr
\noalign{\medskip}
\psi_{t,R}(0)=0\,,\quad
\lim\limits_{|\xi'|\uparrow R}\psi_{t,R}(\xi')=+\infty\,.&\cr}
\ee
It can be proved (see e.g. \cite{LL} for a  more general result in the context of ergodic
problems) that there exists a unique constant $C=C_R$ such that problem \rife{erg}
admits a solution $\psi_{t,R}$, which is also unique. Note that $C_R>0$; moreover, by  a
simple scaling argument, we have
\be\label{scali}
C_R=R^{-\frac{q}{q-1}}C_1\,,\quad \psi_{t,R}=
R^{-\frac{2-q}{q-1}}\psi_{t,1}\left(\frac{|\xi'|}R\right)\,,
\ee
where $C_1$, $\psi_{t,1}$ are the solutions of the same problem in the unit ball
$B_1^{N-1}$. Clearly, we also have that $\psi_{t,R}$ achieves its minimum in zero, hence
$\psi_{t,R}\geq 0$. Consider also $\vfi_{t,L,R}$ solution of 
$$
\cases{
-\vfi ''+  t |\vfi '|^q=\frac{\sqrt{1-t^2}}t
C_R\quad
\hbox{in
$(0,L)$},&\cr
\vfi (0)= \frac Mt\,\quad \vfi (L)= \frac1t(l+C L^{-\frac{2-q}{q-1}}) \,.&\cr}
$$
As in the above case (i), using the concavity of the function $s^{\frac q2}$, one can
check that the function $\bar z=t\vfi_{t,L,R}(\xi_1)+\sqrt{1-t^2}\,\psi_{t,R}(\xi')$ is a
supersolution of \rife{gen} in the cylinder $(0,L)\times B_R^{N-1}$. Moreover, due to 
\rife{lest} and to the properties of $\psi_{t,R}$, we have $\bar z\geq z$ on the
boundary, so that we deduce 
$$
z(\xi)\leq t\vfi_{t,L,R}(\xi_1)+\sqrt{1-t^2}\,\psi_{t,R}(\xi')\,\quad \forall (\xi_1,\xi')\in
(0,L)\times B_R^{N-1}.
$$
In particular for $\xi'=0$ we have $z(\xi_1,0)\leq t\vfi_{t,L,R}(\xi_1)$. Of course we can translate the origin in the $\xi'$--axis, so that we have in fact 
$$
z(\xi)\leq t\vfi_{t,L,R}(\xi_1)\qquad \forall \xi\in \rnp. 
$$
Now let $R$ go to infinity; using \rife{scali} we have that $C_R$ tends to zero, hence we
get
\be\label{upp}
z(\xi )\leq t\vfi_{t,L}(\xi_1)
\ee
where $\vfi_{t,L}$ solves
$$
\cases{
-\vfi_{t,L}''+ t  |\vfi_{t,L}'|^q=0\quad \hbox{in
$(0,L)$},&\cr
\vfi_{t,L}(0)= \frac Mt\,\quad \vfi_{t,L}(L)=\frac1t\,( l+C L^{-\frac{2-q}{q-1}}) \,.&\cr}
$$
As $L$ goes to infinity, $\vfi_{t,L}$ converges to the solution of 
$$
-\vfi_{t}''+ t   |\vfi_{t}'|^q=0\quad \hbox{in
$(0,\infty)$},\quad 
\vfi_{t}(0)= \frac Mt\,\quad \lim\limits_{\xi_1\to +\infty}\vfi_{t}(\xi_1)=\frac 1t\,\min\{l,M\}\,.
$$
Then, inequality \rife{upp} implies, after taking
the limit in $L$, that $z(\xi)\leq t \vfi_t(\xi_1)$ for any $t\in (0,1)$. Note that, in particular, this gives $z\leq M$ on  the whole half space $\rnp$; by definition of $l$, this implies that $l\leq M$.
Now, as    $t$ tends to
$1$, clearly $\vfi_t$ converges to the function $\omega_l(\xi_1)$ defined above. We conclude  that 
\be\label{hai}
z(\xi) \leq \omega_l(\xi_1).
\ee
In order to establish the reverse inequality, let $a\geq 0$, and consider the radial
solutions
$\omega=\omega_{a,R,S}$ of the problems
\be\label{ao}
\cases{-\Delta\omega+ \,|\nabla \omega|^q =0&in $B_{R+S}(0)\setminus
B_R(0)$,\cr
\lim\limits_{\rho\downarrow R} \omega=M\,, \qquad \omega
(&$\!\!\!\!\!\!\!$ R+S)=a\,.\cr}
\ee
Let as before $\eta_R=(-R,\xi')$. We have that the sequence
$\{\omega_{a,R,S}(\xi-\eta_R)\}_R$ is increasing and converges to  a
one--dimensional function $\omega_{a,S}(\xi_1)$ which is the unique solution of
$\omega_{a,S}''=|\omega_{a,S}'|^q$ satisfying $\omega_{a,S}(0)=M$ and $\omega_{a,S}(S)=a$.
As $S$ goes to infinity, we have that $\omega_{a,S}$ converges to $\omega_{a}(\xi_1)$,
which is the unique solution of
$$
\omega_a''=|\omega_a '|^q\quad 
\hbox{in $(0,+\infty)$},\quad \omega_a(0)= M\,,\qquad \lim\limits_{\xi_1\to
+\infty}\omega_a(\xi_1)=a.
$$
In particular, if we know that $z(\xi)\geq a$ for every $\xi\in \rnp$, by comparison we
deduce that $z(\xi)\geq \omega_{a,R}(\xi-\eta_R)$, and then, after letting $R$ and $S$ go
to infinity, that $z(\xi)\geq \omega_a(\xi_1)$. Thus we have the implication
\be\label{prin}
\hbox{$z(\xi)\geq a$ for every $\xi\in \rnp$ implies $z(\xi)\geq \omega_a(\xi_1)$.}
\ee
As a  first step, since $z\geq 0$, this implies that $z\geq \omega_0(\xi_1)$, which
together with \rife{lest} implies
$$
z(\xi)\geq a_1\,:\,= \min
\left[\max\{\omega_0(\xi_1)\,,\,l-C\xi_1^{-\frac{2-q}{q-1}}\}\right]
$$
Note that $0<a_1<l$; applying  \rife{prin} we deduce that $z(\xi)\geq \omega_{a_1}(\xi_1)$
and in particular
$$
z(\xi)\geq a_2\,:\,= \min
\left[\max\{\omega_{a_1}(\xi_1)\,,\,l-C\xi_1^{-\frac{2-q}{q-1}}\}\right]
$$
Iterating this process we define a  sequence of positive real numbers $\{a_n\}$ and a
sequence of functions $\{\omega_{a_n}(\xi_1)\}$ such that  
$$
z\geq \omega_{a_n}(\xi_1)\,,\quad a_n= \min
\left[\max\{\omega_{a_{n-1}}(\xi_1)\,,\,l-C\xi_1^{-\frac{2-q}{q-1}}\}\right]\,.
$$
As $n$ goes to infinity, clearly we have that $a_n\uparrow l$ and $\omega_{a_n}(\xi_1)$
converges to
$\omega_l(\xi_1)$, which allows to conclude that
$$
z\geq \omega_l(\xi_1)\,.
$$
Together with \rife{hai} this concludes the proof.

The case $q=2$ is much simpler. Indeed, if $M<\infty$ it should be noted that the only nonnegative solution of $\omega''= |\omega'|^2$ is the constant $\omega\equiv M$. In particular, one can define $\vfi_{t,L,R}$ as above except for  requiring $\vfi_{t,L,R}(L)=+\infty$; in the limit (in $R$, $L$, $t$ subsequently) one finds that $z\leq M$, while from below one has that   $\omega_{0,S}$ (defined in \rife{ao} for $a=0$) also converges to the constant $M$, so that one gets $z\geq M$, and then    $z\equiv M$. If $M=+\infty$, the function 
$v=e^{-z}$ turns out to be harmonic in $\rnp$ with $v=0$ on $\{\xi_1=0\}$; but $v$ is also asked to satisfy  $0<v\leq 1$, and such a  function cannot exist. 
\qed


\begin{thebibliography}{99}


\bibitem{BE} C. Bandle,    M. Essen, {\it On the solutions of quasilinear elliptic
problems with boundary blow--up}, Symposia Matematica {\bf 35} (1994), 93--111.
     
\bibitem{BG} C. Bandle, E. Giarrusso, {\it Boundary blow--up for semilinear elliptic
equations with nonlinear gradient terms}, Adv. Diff. Equat. {\bf 1} (1996), 133--150.

\bibitem{BaM} C. Bandle, M. Marcus, {\it Large solutions of semilinear elliptic equations:
existence, uniqueness and asymptotic behaviour}, J. Anal. Math. {\bf 58} (1992), 9--24.

\bibitem{BaM2} C. Bandle, M. Marcus, {\it Asymptotic behaviour of solutions and their
derivatives, for semilinear elliptic problems with blowup on the boundary},  Ann. Inst. H.
PoincarŽ Anal. Non LinŽaire {\bf 12}  (1995), no. 2, 155--171. 

\bibitem{D-L} G. Diaz, R.  Letelier, {\it 
Local estimates: uniqueness of solutions to some nonlinear elliptic equations}, 
Rev. Real Acad. Cienc. Exact. F's. Natur. Madrid {\bf 88} (1994), n. 2-3, 171--186.

\bibitem{D-L2} G. Diaz, R.  Letelier, {\it  Explosive solutions of quasilinear elliptic
equations: existence and uniqueness},  Nonlinear Anal. {\bf 20} (1993), n. 2, 97--125.



\bibitem{GNR} M. Ghergu,C. Niculescu, V. Radulescu,  
{\it Explosive solutions of elliptic equations with absorption and non-linear gradient
term}, 
Proc. Indian Acad. Sci. Math. Sci. {\bf 112} (2002), no. 3, 441--451.

\bibitem{G} E. Giarrusso,  {\it Asymptotic behaviour of large solutions of an elliptic
quasilinear equation in a borderline case}, C. R. Acad. Sci. Paris SŽr. I Math. 331
(2000), no. 10, 777--782. 

\bibitem{GT} Gilbarg, N. Trudinger, {\sl Partial Differential Equations of Second Order},
2nd ed., Springer--Verlag, Berlin/New-York, 1983.


\bibitem{K} J.B. Keller, {\it On solutions of $\Delta u=f(u)$}, 
Commun. Pure Appl. Math. {\bf 10} (1957), 503--510.



\bibitem{LL} J.-M. Lasry, P.-L. Lions, {\it  Nonlinear elliptic equations with singular
boundary conditions and stochastic control with state constraints. I. The model problem}, 
Math. Ann. {\bf 283} (1989), n. 4, 583--630. 
 
 
 \bibitem{LMK} A.C. Lazer, P.J. McKenna, {\it Asymptotic behaviour  of solutions of
boundary blow up problems}, Diff. Int. Equ. {\bf 7} (1994), 1001--1019.

\bibitem{Li} P.L. Lions, {\it Quelques remarques sur les problemes elliptiques quasilineaires du second ordre},
 J. Analyse Math. 45 (1985), 234--254.
 
\bibitem{LN} Lowner, L. Nirenberg, {\it Partial differential equations invariant under
conformal or projective transformations.} Contributions to
   analysis (a collection of papers dedicated to Lipman Bers), pp. 245--272, Academic
Press, New York, 1974.


\bibitem{MV} M. Marcus, L. Veron, {\it Uniqueness and asymptotic behaviour of solutions
with boundary blow-up for a class of nonlinear elliptic equations}, Ann. Inst. H.
Poincar\'e {\bf 14} (1997), 237--274.

\bibitem{MV1} M. Marcus, L. Veron, {\it Existence and uniqueness results for large
solutions of general nonlinear elliptic equations}, J. Evolution Equ. {\bf 3}
(2004), 637--652.

 
\bibitem{O} R. Osserman, {\it On the inequality $\Delta u\geq f(u)$}, 
Pacific J. Math. {\bf 7} (1957), 1641--1647.

\bibitem{Ploc} A. Porretta, {\it Local estimates and large solutions for some elliptic
equations with absorption},   Adv. in Diff. Equ.  {\bf 9 }, n. 3/4 (2004), 
329--351.

\bibitem{PoVe} A. Porretta, L. Veron, {\it Symmetry properties of solutions of semilinear 
elliptic equations in the plane}, Manuscripta Math., {\bf 115 } (2004). 239--258.

\bibitem{Ve} L. Veron, {\it Semilinear elliptic equations with uniform blow-up on the
boundary}, J. Anal. Math. {\bf 59} (1992), 231--250.

 \end{thebibliography}
\end{document}